\documentclass{amsart}
\usepackage[latin1]{inputenc}
\usepackage{amsmath}
\usepackage{amsfonts}
\usepackage{amssymb}
\usepackage{amsthm}
\usepackage{calrsfs}
\usepackage{color}
\usepackage{enumitem}
\usepackage{tikz} % Diagram
\usepackage{float}

\usepackage{graphicx} % Graphic
\usepackage{longtable} %Long Table
\usepackage{array} %table length
\usepackage{multirow}

\newcommand{\Z}{\mathbb{Z}}
\renewcommand{\S}{\mathcal{S}}
\newcommand{\N}{\mathbb{N}}

\DeclareMathOperator{\sign}{sign}
\DeclareMathOperator{\Pot}{Pot}

\newtheorem{theorem}{Theorem}
\newtheorem{proposition}[theorem]{Proposition}
\newtheorem{lemma}[theorem]{Lemma}
\newtheorem{cor}[theorem]{Corollary}
\newtheorem{conjecture}[theorem]{Conjecture}

\theoremstyle{remark}

\theoremstyle{definition}
\newtheorem{definition}[theorem]{Definition}

\usepackage{soul}

\title{Rational growth in Torus bundle groups of odd trace}
\author{Seongjun Choi, Meng-Che ``Turbo'' Ho, and Mark Pengitore}
\thanks{The first author is partially supported by NSF grant DMS-1600857 ``Von Neumann Algebras: Rigidity, Applications to Measurable Dynamics, and Model Theory''.}
\thanks{The second author acknowledges support from U.S. National Science Foundation grants DMS 1107452, 1107263, 1107367 ``RNMS: Geometric Structures and Representation Varieties'' (the GEAR Network).}
\date{\today}
\begin{document}

\maketitle

\section{Introduction}

The Cayley graph $\Gamma(G;S)$ of a finitely generated group $G$ with a finite generating subset $S$, when endowed with the word metric $d_S$, allows us to study the group $G$ using both combinatorial and geometric methods. One invariant of the group that arises in this viewpoint is the growth of the size of a ball of radius $n$ in $\Gamma(G;S)$ called the \textbf{growth function} $V(n;S)$ which was introduced by Schwarz and Milnor independently. The motivation for this invariant is geometric as seen in case when $G$ is the fundamental group of a Riemannian manifold. Here, we have that the growth function of $G$ gives a discrete approximation of the volume growth of the universal cover of the manifold \cite{Svr55}. Moreover, it has been shown that the growth function has an exponential growth when the manifold is negatively curved and has polynomial growth when $G$ is virtually nilpotent \cite{Mil68a}\cite{Wlf68}.

%Given two groups $G$ and $H$ with generating subsets $S$ and $X$, respectively, we say two groups have the same \textbf{growth rate} if $V_{G,S}(n/C) \le V_{H,X}(n) \le V_{G,S}(Cn)$ for some constant $C$ from which it easily seen that the growth function is a quasi-isometric invariant of the group. Moreover, there is a rich literature relating the properties of the growth of a group to the structure of the group. (For instance, see \cite{Nek18}.)

Milnor asked whether the growth function is always an exponential function or a polynomial function and whether one can classify every group whose growth exponent $\lim{\log V(n)}/{\log n}$ exists \cite{Mil68b}. It was shown by Tits \cite{tits1972free} that there is no linear group of intermediate growth whereas it is known to be not true in general by work of Grigorchuk \cite{grigorchuk1983milnor}. For the second problem, Bass showed that virtually nilpotent groups always have an integer degree of polynomial growth \cite{bass1972degree}, and Gromov showed that a group can have with polynomial growth if and only if it is virtually nilpotent \cite{Gr81} giving a complete answer to Milnor's question.

For nilpotent groups, one can ask how precisely the growth function $V(n;S)$ behaves like a polynomial over $\mathbb{Q}$. Pansu showed that the limit $\alpha=\lim_{n\rightarrow \infty }V(n;S)/n^d$ exists for nilpotent groups \cite{pansu1983croissance}. Thus, one may study pursue whether $\alpha$ is a rational number. One approach to this question is to study the \textbf{growth series} $\sum V(n;S)t^n$ associated to the growth function and use its analytic properties. If the growth series is a rational function, for instance, one can show that the limit $\alpha$ is an algebraic integer as the coefficients are linearly recursive.

It turns out the rationality of the growth series have stronger implications to the computability of geodesics with respect to a finite generating set. As a corollary to the fact that $V(n)$ is linearly recursive if the growth series is rational, one can see the group $G$ always has a solvable word problem. Cannon showed the growth series for any hyperbolic group is rational by essentially showing that all hyperbolic groups are strongly geodesically automatic \cite{Ca84}. In fact, Neumann and Shapiro showed that the growth series with respect to $S$ is rational when the full language of geodesics in $G$ is regular \cite{neumann1995automatic}. 

Little is known about what groups have rational growth with respect to some finite generating subset. In fact, there are only a handful of classes of finitely generated groups that are known to have rational growth with respect to any finite generating subset.  Virtually abelian groups were shown by Benson \cite{Be83} to have rational growth. For Coxeter groups, it was shown by Paris \cite{paris1991growth}, and for solvable Baumslag-Solitar groups, it was shown by Brazil, and independently by Collins, Edjvet, and Gill \cite{Br94,Co94}. Among these classes of finitely generated groups, hyperbolic groups and virtually abelian groups are also shown to have rational growth with respect to not only some fixed finite generating set, but for all finite generating sets \cite{Ca84,Be83}. This property is known as \textbf{panrationality} and is known to be quite rare. Another known example of a panrational group is the integral Heisenberg group of dimension 3 \cite{Du19}. In general, the rationality of the growth of a group is dependent on the choice of generating subset which was shown by Stoll \cite{St96} in the case of the integral Heisenberg group of dimension 5.  Stoll constructed two generating subsets for this group, one for which the generating function associated to the growth function is rational and the other where the generating function is not rational. Thus, a natural direction of research is to investigate when a finitely generated group of interest admits finite generating subsets that are rational, and we may consider the finite generating sets for which the group has rational growth as algorithmically nicer than other generating sets.

Our interests are in the study of rationality of finitely presented solvable groups. Since Khalarmpovich \cite{Kh81} constructed examples of finitely presented solvable groups of derived length $3$ that have unsolvable word problem, we will focus on metabelian, i.e.\ solvable of length 2, groups in this paper. In particular, we will focus on a family of torus bundle groups given by $G = \Z^2 \rtimes \Z$. Parry showed that when the trace of the action is even and at least 4, these groups have rational growth with respect to some generating set \cite{Pa07}. We extend Parry's result to the traces that are odd and at least 5.

\begin{theorem}
Let $A = \begin{bmatrix} 0 & -1 \\ 1 & 2k+1\end{bmatrix}\in SL(2,\Z)$ where $2k+1 \ge 5$. Then $\Z^2\rtimes_A\Z$ has rational growth with respect to the standard generating set.
\end{theorem}

It was shown that certain groups of the form $\Z^m\rtimes \Z$ have a finite-index subgroup which has rational growth in some finite generating set \cite{Pu06,Br19}. Our work improves on this result for $m = 2$ as we do not require passing to a finite-index subgroup. Moreover, we are able to generalize and to streamline the proof given by Parry. In particular, we define an invariant called the \textbf{potential} that makes the proof cleaner. We conjecture that our argument can be applied to all torus bundle groups to get the following:

\begin{conjecture}
Any torus bundle group $G = \Z^2 \rtimes_A \Z $ where $A \in SL(2,\Z)$ has rational growth with respect to some finite generating set.
\end{conjecture}

%\opinion{Add in description of the cases for $A = \begin{bmatrix} 0 & -1 \\ 1 & d\end{bmatrix}$: $d = 0$ is abelian, $d = 1$ is nilpotent, $d = 2, 3$ are degenerate, and negative $d$ via a bijection.}

The general structure of this work is as follows: in section 2, we first recall the definition of a torus bundle group and its structure, giving a description of the word length with a specific choice of generators. In section 3, we by following \cite{Pa07} reduce the problem of finding the word with the shortest length representing a group element into an optimization problem in the Laurent polynomial ring over the integers. Moreover, we build a correspondence between a group element $\gamma$ and a Laurent polynomial $P$ called a $n$-minimal polynomial. In order to quantify the $n$-minimal polynomials, we first define a canonical form called the $n$-reduced polynomial using the potential. In section 4, we prove that there is a deterministic way of rewriting a given ordinary polynomial into the desired form. Using this machinery, we show that any $n$-minimal polynomial can be split into two $n$-reduced polynomials. In order to count all $n$-reduced polynomials, we define a successor map in section 5 that gives a natural ordering on the set of $n$-reduced polynomials and partitions the behavior of the successor function. This allows us to build a recursive relation which will be encoded as a transition matrix between partitions, giving the growth function.

%\section{Definitions}

%A group $G$ is said to have a \textbf{rational growth} if the growth series $B(t)=\sum_{n=0} V_G(n) t^n$ is a rational function with respect to some generating subset. In this paper, we investigate the rationality of a family of groups called the \textbf{torus bundle groups} which was first investigated by W.~Parry \cite{Pa07}.

\section{Torus bundle groups}
In this section, we briefly give a definition of the class of torus bundle groups we are interested in and a description of their word metrics with respect to a specific choice of finite generating set. Let $T$ be the automorphism of $\Z^2$ given by the matrix
$$
T=\begin{bmatrix}
0 & -1 \\
1 & d
\end{bmatrix}
$$
with respect to the standard basis
$
a = [\begin{smallmatrix}
1 \\ 0
\end{smallmatrix}], \:
b = [\begin{smallmatrix}
0 \\ 1
\end{smallmatrix}]
$
of $\Z^2$. Let $\langle t \rangle$ be the infinite cyclic group generated by $t$. We define the action of $t$ on $\Z^2$ by $
t \: x \: t^{-1} = T \cdot x$ for any $x \in \mathbb{Z}^2$ and form the semidirect product $G = \Z^2 \rtimes \langle t \rangle$. We will view elements of $G$ as ordered pairs $\left(x,t^n\right)$ where $x \in \Z^2$ and $n \in \Z$. Thus, we may write $G$ as
$$
G = \langle a,b,t \mid [a,b]=1, a^t = b, b^t = a^{-1}b^{d} \rangle$$
where $g^t = t \: g \: t^{-1}$. In particular, we will investigate the rational growth of $G$ with respect to the generating set $S = \{a,b,t\}$.

\subsection{Length of elements in $G$}
Let $g=(x,t^n) \in G$. Through a rewriting process where we collect consecutive factors which are either $a$'s or $b$'s but not collecting consecutive factors which are either $t$'s or their inverses, we may write
$$
g = t^{n_0} x_1 \: t^{n_1} \: x_2 \: t^{n_2} \: \cdots x_{k-1} t^{n_{k-1}} \: x_k \: t^{n_k}
$$
where $h \in \mathbb{N}$, $n_0,n_k \in \{-1,0,1\}$, $n_1,\dotsc, n_{k-1} \in \{-1,1\}$, and $x_1,\dotsc,x_k \in \Z^2$. We may obtain the word length of $g$ with respect to $S$, denoted as $\|g\|$, by minimizing the sum
\begin{equation}\label{geodesic_length_1}
\sum_{i=0}^{k}|n_i| + \sum_{i=1}^k |x_i|,
\end{equation}
where $|x_i|$ is the length of $x_i$ in $\Z^2$ with respect to the generating subset $\{a,b\}$. We rewrite the given expression for $g$ as
\begin{equation}\label{element_expression_1}
g = (t^{n_0} \: x_1 \: t^{-n_0}) \cdot (t^{n_0 + n_1} \: x_2 \: x^{-n_0 - n_1}) \cdots (t^{\sum_{i=0}^{k-1} n_i} \: x_k \: t^{-\sum_{i=0}^{k-1} n_i}) \cdot t^{\sum_{i=0}^{k}n_i}.
\end{equation}
One can observe that each of these terms lie in $\Z^2$ and that
$
t^{m} \: x \: t^{-m} = T^m \cdot x
$
for all $m \in \Z$ and $x \in \Z^2$. Therefore, we have that
$$
(t^{n_0} \: x_1 \: t^{-n_0}) \cdot (t^{n_0 + n_1} \: x_2 \: x^{-n_0 - n_1}) \cdots (t^{\sum_{i=0}^{k-1} n_i} \: x_k \: t^{-\sum_{i=0}^{k-1} n_i}) \in \Z^2.
$$
By the above construction, we have that $n = \sum_{i=0}^k n_i$. Let 
$$
q = \text{max}\{0,n_0,n_0 + n_1, n_0 + n_1 + n_2,\dotsc, n\}
$$
and 
$$
p = \text{max}\{0,-n_0,-n_0 - n_1, -n_0 - n_1 - n_2,\dotsc, -n\}.
$$

Assume that $n \geq 0$. By combining terms according to the partial sums $n_0 + \cdots + n_i$, we will not increase the length of the expression in Equation \ref{geodesic_length_1}. Thus, we may rewrite Equation \ref{element_expression_1} as
\begin{eqnarray*}
g & = & \left(t^{-p} \: y_{-p} \: t^p \right) \cdot \left( t^{-p + 1} \: y_{-p+1} \: t^{p-1}\right) \cdots \left( t^q \: y_q t^{-q} \right) \cdot t^n \\
& = &  \left( \sum_{i=-p}^q T^{i}(y_i),t^n \right)
\end{eqnarray*}
for elements $y_{-p},\dotsc, y_q \in \Z^2$. We furthermore obtain $\|g\|$ by minimizing 
$$
2p + 2 q - n + \sum_{i=-p}^q |y_i|.
$$

We may suppose that $y_i = r_i \: a + s_i \: b$ with $r_i,s_i \in \Z$. Given that $T(a) = b$, if $i > 0$, we may write $T^i(y_i)$ as
$$
T^i(y_i) = r_i \: T^{i}(a) + s_i T^{i}(b) = r_i \: T^{i-1}(b) + s_i \: T^{i}(b).
$$
Additionally, if $i < 0$, then we may write $T^i(y_i)$ as
$$
T^i(y_i) = r_i \: T^i(a) + s_i T^i(b) = r_i T^{i}(a) + s_i T^{i+1}(a).
$$
We have a similar situation for $n \leq 0$. Thus, if $g = (x,t^n) \in G$, then $\|g\|$ is the minimal value of
$$
2p + 2q - |n| + \sum_{i=0}^p |r_i| + \sum_{j=0}^{q}|s_j|
$$
over all representatives of $x$ given by
$$
x = \sum_{i=0}^p r_i T^{-i}(a) + \sum_{j=0}^q s_j T^{j}(b)
$$
where $r_0,\dotsc, r_p,s_0,\dotsc,s_q$ are integers, $p \geq \text{max}\{0,-n\}$, and $q \geq \text{max}\{0,n\}$. Given that $a = T^{-1}(b)$, we may instead write $\|g\|$ as the minimal value of
$$
2p + 2q - |n| + \sum_{j=-p}^q |c_j|
$$
over all representations of $x$ given by
$$
x = \sum_{j=-p-1}^q c_j T^j(b),
$$
where $p \geq \text{max}\{0,-n\}$, $q \geq \text{max}\{0,n\}$, and $c_{-p-1},\dotsc, c_q$ are integers.

\section{Polynomial representatives}
%Comment [1]

Let $F(X) = \sum_{j=-p-1}^q c_j X^j$ be a Laurent polynomial in $\Z[X, X^{-1}]$. We refer to $\sum_{j=0}^q c_j X^j$ as the \textbf{polynomial part} of $F(X)$ and $\sum_{j=-p-1}^{-1} c_j X^j$ as the \textbf{principal part} of $F(X)$. Let $(x, t^n) \in G$. Then $F(X)$ will be called a \textbf{representative of $x$} if $x = F(T)(b)$ and will be said to \textbf{represent $x$}. By eliminating and introducing terms with $c_j = 0$,  we may assume that $p \geq \text{max}\{0, -n\}$ with $c_{-p-1} \neq 0$ if $p > \text{max}\{0, - n\}$ and $q \geq \text{max}\{0, n\}$ with $c_q \neq 0$ if $q > \text{max}\{0, n\}$. It is straightforward to see that the Laurent polynomial $F(X)$ and integer $n$ uniquely determine such integers $p$ and $q$. In this situation, the non-negative integer
$$
L_n(F(X)) = 2p + 2q - |n| + \sum_{j = -p -1}^q |c_j|
$$
will be called the \textbf{$n$-length} of $F(X)$. Furthermore, if $F(X)$ represents $x$, then $F(X)$ will be called an \textbf{$n$-minimal representative} of $x$ if $L_n(F(X)) = \left\|(x,t^n)\right\|$. Since $L_n(F(X))$ is a positive-valued function on the integers unless $F(X)=0$, we have the following theorem:

\begin{theorem}\label{n-min theorem}
Suppose that $g=(x,t^n)$ is an element of $G$. Then $g$ has a $n$-minimal representative.
\end{theorem}

We note the relation $b^t = (b^{-1})^{t^{-1}}b^{2k+1}$ in this context becomes $T-(2k+1)I+T^{-1}=0$ and thus translates to minimizing $L_n$ over all polynomials representatives $F(X)$ up to adding multiples of $(X^2-(2k+1)X+1)$. One can easily observe that when given a representative $F(X)$ that has a large coefficient, one may cancel with the relation $(X^2-(2k+1)X+1)$ and yield a smaller coefficient. For example, the polynomial $X^2+(k+3)X+3$ is not preferable to another polynomial representative $2X^2+(-k+2)X+4$ which is obtained by adding $(X^2-(2k+1)X+1)$. In fact, one can use this cancellation repeatedly to deduce that large coefficients cannot be adjacent to each other too often in a $n$-minimal polynomial. We will later call this phenomenon the \textbf{global rules}. For example, consider $$X^6+(k+1)X^5+(k-1)X^4+kX^3+kX^2+kX+2.$$
By adding $$X^6+(-2k)X^5+(-2k+1)X^4+(-2k+1)X^3+(-2k+1)X^2+(-2k)X+1,$$ we get the polynomial $$2X^6+(-k+1)X^5+(-k)X^4+(-k+1)X^3+(-k+1)X^2+(-k)X+3$$ which has a shorter $n$-length. More precise conditions on the coefficients will be described below.

\subsection{$n$-reduced polynomials}
In this subsection, we define \textbf{$n$-reduced polynomials} in $\mathbb{Z}[X]$. Later, we show in {Theorem \ref{decomposition.thm}} that any $n$-minimal Laurent polynomial can be written as the sum of two $n$-reduced polynomials, namely, the polynomial part and the principal part. Since the set of $n$-reduced polynomials has a natural grading given by the degree of a polynomial, we will be using $n$-reduced polynomials instead of $n$-minimal Laurent polynomials for counting purposes. For convenience, we denote a polynomial $F(X)=\sum_{i=0}^{m}c_i X^i$ by a string of coefficients $(c_m,c_{m-1},\dotsc,c_{1},c_0)$ called the \textbf{word} where indices are understood as the degree of $X$. A contiguous substring of this will be called the \textbf{subword} of a word. We adopt the following convention of $\sign$ functions:

\begin{center}
\begin{tabular}{ c c c }
$
\sign(x) = \begin{cases}
1 & x > 0 \\
0 & x = 0 \\
-1 & x <0
\end{cases},
$
&
$\sign^+(x) = \begin{cases}
1 & x \geq 0 \\
-1 & x < 0 
\end{cases},$
&
$\sign^-(x) = -\sign^+(-x)$
\end{tabular}
\end{center}
with the inequality $\sign^-(x+1)\leq \sign^+(x)$ for all $x\in \mathbb{Z}$. We start with the following definition.

\begin{definition}
For a word $(c_{m},\dotsc, c_{m-i})$, we define its \textbf{potential} as $$\Pot(c_m,\dotsc, c_{m-i}) = \sum_{c_j\in(c_m,\dotsc, c_{m-i})} (2k-1)-2|{c_j}|.$$
In particular, when the string only contains coefficients whose absolute values are $k-2$, $k-1$, $k$ or $k+1$, then the potential can be written as

$$\Pot(c_m,\dotsc, c_{m-i})=3N_{k-2}+N_{k-1}-N_{k}-3N_{k+1}$$
where $N_t = |\{ j \: | \: c_{j} = |t| \} |$.
\end{definition}

By definition, the potential is additive under concatenating two strings. In order to describe the change in the potential under any rewriting, we analyze the potential change in different cases on a single coefficient $a$ with $|a|<k+2$:

\begin{itemize}
    \item $\Pot(a+1)=\Pot(a)-2\sign^{+}(a)$
    \item $\Pot(a-1)=\Pot(a)+2\sign^{-}(a)$
    \item $\Pot(a-(2k-1))=-\Pot(a) \quad (a>0)$
    \item $\Pot(a+(2k-1))=-\Pot(a) \quad (a<0)$
    \item $\Pot(a-2k)=-\Pot(a)-2 \quad (a>0)$
    \item $\Pot(a+2k)=-\Pot(a)-2 \quad (a<0)$
\end{itemize}
which clearly follows from definition.

The potential serves as a test function that determines whether a rewriting is needed. We are now ready to define $n$-reduced polynomials.

\begin{definition}
Let $n$ be a fixed integer and $P(X)=(c_m,\dotsc,c_1,c_0)$ be a polynomial of degree $m$. $P(X)$ is \textbf{$n$-reduced} if the following rules are satisfied.
\end{definition}

\begin{center}
    \begin{longtable}{|c|}\hline
    {\parbox{0.9\textwidth}{\vspace{0.5em}
        \textbf{Rule 1} (Local top rule) For $m\ge n$, \begin{align*}
            |c_m|&\leq k+2,\\
            |c_m|&\leq k+1 \quad \text{ if } \sign(c_m \cdot c_{m-1})< 0.
        \end{align*}\vspace{0.5em}
        }
    }\\\hline
    {\parbox{0.9\textwidth}{\vspace{0.5em}
        \textbf{Rule 2} (Local non-top rule) For $i=m< n$ or $i<m$,
        \begin{align*}
            |c_i|&\leq k+1, \\
            |c_i| &\leq k &\text{ if } \sign(c_{i+1} \cdot c_{i})< 0\text{  or }\sign(c_i \cdot c_{i-1})< 0,\\
            |c_i|&\leq k-1 &\text{ if } \sign(c_{i+1} \cdot c_{i})< 0\text{ and }\sign(c_i \cdot c_{i-1})< 0.
        \end{align*}\vspace{0.5em}
        }
    }\\\hline
    {\parbox{0.9\textwidth}{\vspace{0.5em}
        \textbf{Rule 3} (Local top rule) For $m> n$, \begin{align*}
            (c_m,c_{m-1})&\neq \pm (1,-k), \\
            (c_m,c_{m-1})&\neq \pm (1,-k+1) \quad \text{ if } \sign(c_{m-1} \cdot c_{m-2})< 0.
        \end{align*}\vspace{0.5em}
        }
    }\\\hline
    {\parbox{0.9\textwidth}{\vspace{0.5em}
        \textbf{Rule 4} (Global top rule) \par
        For $m>n$ and $i\ge1$, and $$(c_m,\dotsc, c_{m-i}) = \pm(1,c'_{m-1},\dotsc,c'_{m-i})$$ where $c^\prime_{m-1} \in \{-k + 1, - k +2\}$, $c'_{m-j}\in \{-k+1,-k\}$ for $1<j<i$, and $c'_{m-i} \in \{-k, -k-1\}$.
        $$\Pot(c'_{m-1},\dotsc,c'_{m-i}) \geq 2\quad\text{ if }\sign(c_m \cdot c_{m- i-1}) > 0,$$ $$ \Pot(c'_{m-1},\dotsc,c'_{m-i})> 0\quad\text{ if }\sign(c_m \cdot c_{m -i - 1}) \le 0.$$
        }
    }\\\hline
    {\parbox{0.9\textwidth}{\vspace{0.5em}
        \textbf{Rule 5} (Global top rule) \par
        For $m\ge n$, $i\ge1$, and $$(c_m,\dotsc, c_{m-i}) = \pm(c'_m,c'_{m-1},\dotsc,c'_{m-i})$$ where $c'_m \in \{k+1, k+2\}$, $c'_{m-i'}\in \{k-1,k\}$ for $i' \neq i$, and $c'_{m-i}\in \{k,k+1\}$.
        $$\Pot(c'_{m},\dotsc,c'_{m-i})> -6 \quad\text{ if }\sign(c_m \cdot c_{m- i-1}) \geq 0,$$ $$\Pot(c'_{m},\dotsc,c'_{m-i})\geq -4 \quad\text{ if }\sign(c_m \cdot c_{m -i - 1}) <0.$$
        }
    }\\\hline
    {\parbox{0.9\textwidth}{\vspace{0.5em}
        \textbf{Rule 6} (Global non-top rule) \par
        For $0\leq l<j$ and either $0<j<m$ or $0<j = m<n$, and $$(c_j,c_{j-1},\dotsc,c_{l})=\pm(c'_j,c'_{j-1},\dotsc,c'_{l})$$ where $c'_{j},c'_{l} \in \{k,k+1\}$,\par and $c'_{s} \in \{k-1,k\}$ for $s \neq j,l.$
        $${\Pot(c_j,\dotsc,c_{j-i}) > - 3 - \sign^+(c_j \cdot c_{j+1})\:\:\text{ if }\sign(c_j \cdot c_{l-1}) \geq 0}, $$ $${\Pot(c_j,\dotsc,c_{j-i}) \geq -1 - \sign^+(c_j \cdot c_{j+1}) \:\:\text{ if }\sign(c_j \cdot c_{l-1}) < 0}.$$
        }
    }\\\hline
    \end{longtable}
\end{center}

\textbf{Rule 1}, \textbf{2} and \textbf{3} are called \textbf{local rules} which determine whether a single coefficient requires rewriting. \textbf{rule 4}, \textbf{5} and \textbf{6} are called \textbf{global rules} which determine whether a string of coefficients requires rewriting. The local rules can be considered as the degenerate cases of the corresponding global rules. Note that the local rules are mutually exclusive, and \textbf{rule 4} and \textbf{5} are mutually exclusive. For the global rules, a subword $(c_j,c_{j-1},\dotsc,c_{j-i})$ used in the potential condition will be called the \textbf{subword associated to the rule}. For the local rules, the associated subword will be understood as a single coefficient $\pm(k-1)$, $\pm k$, $\pm(k+1)$, or $\pm(k+2)$ in question. For $F(X)\in \mathbb{Z}[X^{-1}]$, we say $F(X)$ is \textbf{$n$-reduced} if $F(X^{-1})\in \mathbb{Z}[X]$ is $n$-reduced.

For convenience, we introduce an alternative definition of global rules.
\begin{lemma}
Assuming the same conditions for subwords, the potential condition for the global rules can be restated as follows:
\end{lemma}

\begin{itemize}
    \item \textbf{Rule 4}
    $$\Pot(c'_{m-1},\dotsc,c'_{m-i}) \geq 1+\tfrac{1}{2}\sign^-(c_m\cdot c_{m-i-1}).$$
    \item \textbf{Rule 5}
     $$\Pot(c'_{m},\dotsc,c'_{m-i})\geq -5- \tfrac{1}{2}\sign^+(c_m \cdot c_{m-i-1}).$$
    \item \textbf{Rule 6}
    $$\Pot(c_j,c_{j-1},\dotsc,c_{l}) \geq -2- \sign^+(c_j\cdot c_{j+1})- \tfrac{1}{2}\sign^+(c_j \cdot c_{l-1})$$
\end{itemize}

\noindent Since $\Pot(\cdot)$ is integer valued, there is no distinction between strict inequality and inequality. This allows us to reduce the number of cases one needs to check.  We will use this alternative definition when appropriate.

We note that not all $n$-minimal polynomials are $n$-reduced; thus, we may have to rewrite the polynomial. Thus, for each violation, we assign \textbf{rewriting rules} that decreases $n$-length.

\begin{definition}
Suppose that $P(X)=(0,c_m,\dotsc,c_0)$ is not $n$-reduced. For each violation of the rules, we define \textbf{the rewriting associated to the given rule} as the following operations.
\end{definition}

\begin{enumerate}
    \item \textbf{Rule 1}
    $$(0,c_m,c_{m-1})-\sign(c_m)(-1,2k+1,-1).$$
    \item \textbf{Rule 2} $$(c_{i+1},c_i,c_{i-1})-\sign(c_i)(-1,2k+1,-1).$$
    \item \textbf{Rule 3}
    $$(1,c_{m-1},c_{m-2})+(-1,2k+1,-1)\text{, or}$$ $$(-1,c_{m-1},c_{m-2})+(1,-2k-1,1).$$
    \item \textbf{Rule 4}\\
    For the subword associated to the rule $(c_{m-1},\dotsc,c_{m-i})$,
    \begin{align*}
        (1,c'_{m-1},\dotsc,c'_{m-i},c'_{m-i-1})+(-1,2k,2k-1,\dotsc,2k-1,2k,-1).
    \end{align*}
    \item \textbf{Rule 5}
    For the subword associated to the rule $(c_{m},\dotsc,c_{m-i})$,
    \begin{align*}
        (0,c_{m},\dotsc ,c_{m-i},c_{m-i-1})-\sign(c_m)(-1,2k,2k-1,\dotsc,2k-1,2k,-1).
    \end{align*}
    \item \textbf{Rule 6}\\
    For the subword associated to the rule $(c_j,c_{j-1},\dotsc,c_{l})$, 
    \begin{align*}
        (c_{j+1},c_{j}\cdots,c_{l},c_{l-1})-\sign(c_j)(-1,2k,2k-1,\dotsc,2k-1,2k,-1).
    \end{align*}
\end{enumerate}

Since we allow the rewriting rules to take place at the constant term $c_0$, the resulting polynomial may have a non-zero coefficient of $X^{-1}$. For example, if $F(X)=k+2$, the rewriting will give us $F'(X)=X+(-k+1)+X^{-1}$. However, if we start with $F(X)\in\mathbb{Z}[X]$, after rewriting, the coefficient for $X^{-1}$ can only be $0$ or $\pm 1$.

\begin{lemma}\label{length.change}
The rewriting associated to the rules will not increase the $n$-length of the sequence.
\end{lemma}
\begin{proof}
For local rules, we check the $n$-length change for each case.
\begin{itemize}
    \item \textbf{Rule 1} The associated rewriting rule is given by
$$(0,k+3,b,\dotsc) + (1,-2k-1,1) = (1,-k+2,b+1,\dotsc).$$ The $n$-length is changed by $3-5+\sign^+(b) < 0$.
    \item \textbf{Rule 2} The associated rewriting rule is given by
$$(\cdots,b,k+2,c,\dotsc) + (1,-2k-1,1) = (\cdots, b+1,-k+1,c+1,\dotsc).$$ The $n$-length is changed by $\sign^+(b) -3 +\sign^+(c) < 0$.
    \item \textbf{Rule 3} The associated rewriting rule is given by 
$$(1,-k,c,\dotsc) + (-1,2k+1,-1) = (0, k+1,c-1,\dotsc).$$
The $n$-length is changed by $-3 +1 \pm1 < 0$.
\end{itemize}
For global rules, suppose that $c_i<0$ is a coefficient in the subword associated to the rule. We see that adding $(2k-1)$ to $c_i$ for all possible values of $c_i$ yields
    \begin{align*}
        -k-1&\mapsto k-2\\
        -k&\mapsto k-1\\
        -k+1&\mapsto k\\
        -k+2&\mapsto k+1,
    \end{align*}
and for each case, the absolute value of the coefficient changes by $-3$, $-1$, $1$ and $3$ respectively, which is precisely how the potential was defined. The rest of the changes in $n$-length can be argued similarly to those of the local rules.
\end{proof}

\section{Rewriting Procedures}

In this section, we investigate the relationship between $n$-minimal polynomials and $n$-reduced polynomials. As noted, while some rules are mutually exclusive, it is possible that a word violates multiple rules at the same time. Even worse, after rewriting, another rule violation may appear somewhere else. For that reason, it is not immediately clear how any polynomial can be uniquely rewritten to a $n$-reduced polynomial. In this section, we show that there is a unique sequence of rewriting rules determined by the given polynomial that gives a $n$-reduced polynomial.

We begin by defining the basic unit of the procedure.
\begin{definition}
If a word violates a unique rule, it is called a \textbf{poison word}. The shortest subword that violates a unique rule is called the \textbf{minimal poison subword}. A minimal poison subword that has the smallest leading degree is called the \textbf{rightmost minimal poison subword}.
\end{definition}
If a word violates multiple rules simultaneously, then there is a shorter subword in it that violates only one of the rules. By the natural ordering by inclusion, we know there has to be a minimal subword, which includes the case when a single coefficient violates one of the local rules. As noted, the local rules are mutually exclusive. Therefore, the minimal poison subword must exist. Although there may be multiple minimal poison subwords, there can be only one rightmost minimal poison subword.

\begin{definition}
Let $(\cdots,c_i,\dotsc)$ and $(\cdots,c_j,\dotsc)$ be two subwords that violate any of the rules. If there is no common subword that is contained in both, then two subwords are said to be \textbf{disjoint}. In the same manner, two rewriting rules are disjoint if their associated subwords are disjoint.
\end{definition}

Two disjoint rules do not necessarily commute unless their associated subwords are separated by another nontrivial subword. Thus, the ordering of the rules is important.

The main idea is as follows: given a word, we first find the rightmost minimal poison subword to rewrite, and then we choose another rightmost minimal poison subword after rewriting and so forth. This makes it a deterministic process. In order for this proof to work, we have to show that this sequence of rewriting progresses to the left without backtracking.

\begin{theorem}\label{unique.write}
Let $(s_i,\dotsc,s_{i-k})$ be the rightmost minimal poison subword of a given word. Let $(s'_i,\dotsc,s'_{i-k})$ be the subword you obtain from rewriting, and let $(s'_j,\dotsc,s'_{j-l})$ be the rightmost minimal poison subword for the new word. Then $j>i$ and $j-l>i-k$. Subsequently, rewriting rules are eventually disjoint from the initial rewriting.
\end{theorem}

\begin{figure}[ht]
	\includegraphics[scale=0.5]{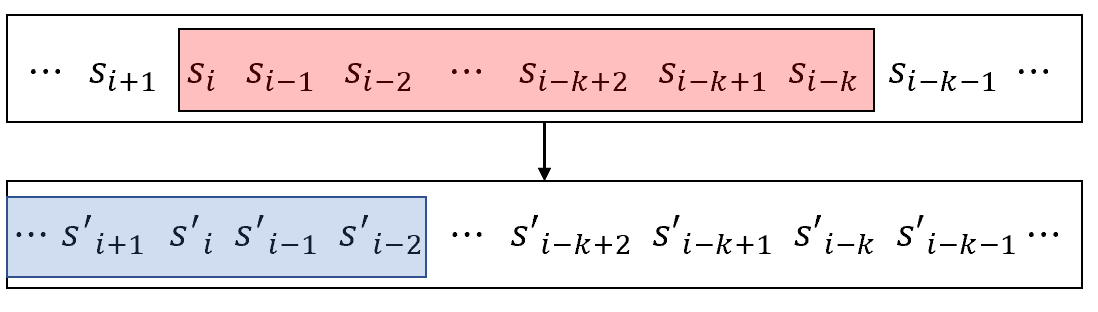}
	\caption{Given a word (top), we perform a rewriting with respect to the minimal poison subword (in red). After rewriting, a new minimal poison subword always appears on the left (in blue).}
\end{figure}

\begin{proof}
We rule out all of the other possible cases using the potential of subwords.
\begin{enumerate}[wide, labelindent=0pt]
\item The subword $(s'_j,\dotsc,s'_{j-l})$ is contained in the subword $(s'_i,\dotsc,s'_{i-k})$.\par
	\begin{figure}[ht]
		\includegraphics[scale=0.5]{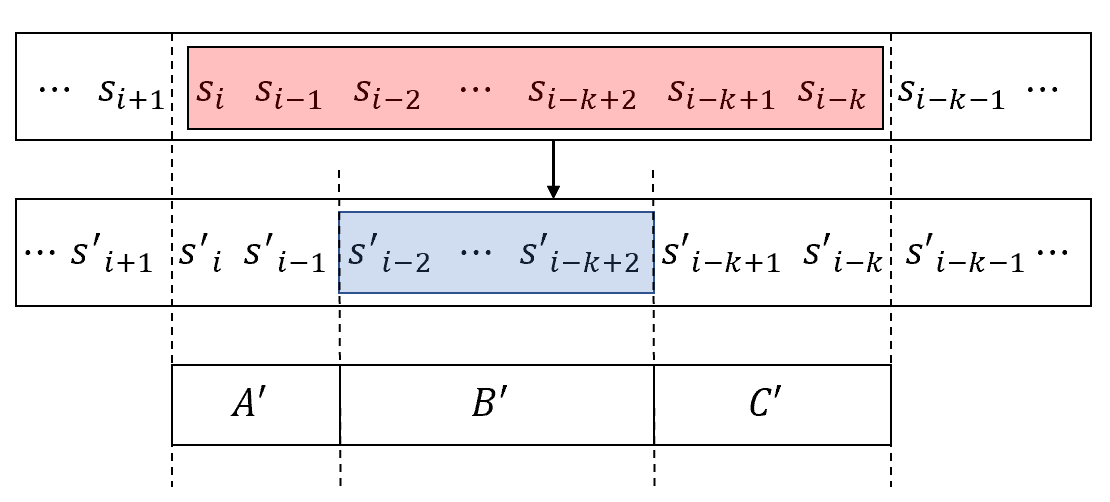}
		\caption{The new minimal poison subword cannot be contained in where the previous minimal poison was.}\label{contained-case}
	\end{figure}
	We first show that the new poison subword $(s'_j,\dotsc,s'_{j-l})$ cannot be contained in the subword $(s'_i,\dotsc,s'_{i-k})$ (see Figure \ref{contained-case}.) For simplicity, denote $(s'_i,\dotsc,s'_{i-k})=(A'|B'|C')$ as 3 disjoint subwords concatenated into one word where $(s'_j,\dotsc,s'_{j-l})=(B')$, and correspondingly, denote $(s_i,\dotsc,s_{i-k})=(A|B|C)$. The rewriting rule in this case is \begin{align*}
	&(s_{i+1}, |A|B|C|, s_{i-k-1})-\sign(s_j)(-1,2k,2k-1,\dotsc,2k-1,2k,-1)\\
	&=(s'_{i+1}, |A'|B'|C'|, s'_{i-k-1}),
\end{align*} and all coefficients in $(A), (B)$ and $(C)$ have the same sign. Moreover, the sign of the coefficients will be denoted as $\sign(A),$ $\sign(B)$, and $\sign(C)$, respectively. Obviously, $(A)$ and $(C)$ cannot be empty at the same time.

We argue in the following way: because $(A'|B'|C')$ was obtained from rewriting, there's a bound on $\Pot(A'|B'|C')$ from $\Pot(A|B|C)$. Since $(B')$ violates one of the rules, we can combine with alternate bounds on $\Pot(A')$ and $\Pot(C')$ to get another bound on $\Pot(A'|B'|C')$ contradicting the previous bound on $\Pot(A'|B'|C')$.

Suppose that both $(A)$ and $(C)$ are nonempty and $(A|B|C)$ violates \textbf{rule 6}. Recall that $(s_i,\dotsc,s_{i-k})$ satisfies \textbf{rule 6} if
\begin{eqnarray*}
	\Pot(s_i,s_{i-1},\dotsc,s_{i-k}) &\geq& -2- \sign^+(s_i \cdot s_{i+1})- \tfrac{1}{2}\sign^+(s_i \cdot s_{i-k-1}).
\end{eqnarray*}
Since $(s_i,\dotsc,s_{i-k})=(A|B|C)$ violates \textbf{rule 6}, we instead have
\begin{align*}
	\Pot(A|B|C) &\leq -2- \sign^+(s_i \cdot s_{i+1})-\tfrac{1}{2}\sign^+(s_i \cdot s_{i-k-1}).
\end{align*}
Since $(A|B|C)$ is minimal, none of the subwords inside $(A|B|C)$ violates any rules, in particular, \textbf{rule 6}. Thus, we have
\begin{align*}
	\Pot(A) &\geq -2- \sign^+(s_i \cdot s_{i+1})- \tfrac{1}{2}\sign^+(A\cdot B),\\
	\Pot(A) &\geq -2- \sign^+(s_i \cdot s_{i+1}),\\
    \Pot(C) &\geq -2- \sign^+(B \cdot C)-\tfrac{1}{2}\sign^+(s_i\cdot s_{i-k-1}),\\
    \Pot(C) &\geq -3-\tfrac{1}{2}\sign^+(s_i\cdot s_{i-k-1}).
\end{align*}
After rewriting, $(B')$ is the rightmost minimal poison subword violating \textbf{rule 6}. Therefore, we have $$\Pot(B') \leq -2- \sign^+(A'\cdot B')- \sign^+(B'\cdot C')=-4.$$
Because $(A'|B'|C')$ was obtained from $(A|B|C)$ by rewriting, we have
\begin{align*}
    \Pot(A')&=-\Pot(A)-2\leq \sign^+(s_i \cdot s_{i+1}),\\
    \Pot(C')&=-\Pot(C)-2\leq 1+\tfrac{1}{2}\sign^+(s_i\cdot s_{i-k-1}),\\
    \Pot(A'|B'|C')&=-\Pot(A|B|C)-4\geq -2+\sign^+(s_i \cdot s_{i+1})+\tfrac{1}{2}\sign^+(s_i \cdot s_{i-k-1}).
\end{align*}
But $\Pot(B')\leq -4$. Hence, $$\Pot(A'|B'|C')\leq -3+\sign^+(s_i \cdot s_{i+1})+\tfrac{1}{2}\sign^+(s_i\cdot s_{i-k-1})$$ which is a contradiction. If $(A)$ is empty, then we have
\begin{align*}
    \Pot(B')&\leq -3-\sign^+(s'_i \cdot s'_{i+1}),\\
    \Pot(C')&\leq 1+\tfrac{1}{2}\sign^+(s_i\cdot s_{i-k-1}),\\
    \Pot(B'|C')&=-\Pot(B|C)-4\geq -2+\sign^+(s_i \cdot s_{i+1})+\tfrac{1}{2}\sign^+(s_i\cdot s_{i-k-1}).
\end{align*}
However, since $\sign^-(s_i \cdot s'_{i+1})\leq \sign^+(s_i \cdot s_{i+1})$, we have that
\begin{align*}
    \Pot(B'|C')&\leq -2-\sign^+(s'_i \cdot s'_{i+1})+\tfrac{1}{2}\sign^+(s_i\cdot s_{i-k-1}),\\&=-2+\sign^-(s_i \cdot s'_{i+1})+\tfrac{1}{2}\sign^+(s_i\cdot s_{i-k-1}),\\
    &\leq-2+\sign^+(s_i \cdot s_{i+1})+\tfrac{1}{2}\sign^+(s_i\cdot s_{i-k-1}),\\
    &\leq \Pot(B'|C'),
\end{align*}
and the equality can only happen when $\sign^+(s_i\cdot s_{i-k-1})=0$ which is impossible. When $(C)$ is empty instead, we have
\begin{align*}
    \Pot(A')&=-\Pot(A)-2\leq \sign^+(s_i \cdot s_{i+1}),\\
    \Pot(B')&\leq -3+\tfrac{1}{2}\sign^+(s'_i\cdot s'_{i-k-1}),\\
    \Pot(A'|B')&=-\Pot(A|B)-4\geq -2+\sign^+(s_i \cdot s_{i+1})+\tfrac{1}{2}\sign^+(s_i\cdot s_{i-k-1}),\\
    \Pot(A'|B')&\leq -3+\sign^+(s_i \cdot s_{i+1})+\tfrac{1}{2}\sign^+(s'_i\cdot s'_{i-k-1}),\\
    &\leq -3+\sign^+(s_i \cdot s_{i+1})+\tfrac{1}{2}\sign^+(s_i\cdot s_{i-k-1})<\Pot(A'|B').
\end{align*}

Now we check for the \textbf{rule 4} and \textbf{rule 5} altogether. As before, $(A)$ and $(C)$ cannot be both empty, as these two rules are mutually exclusive. Suppose that both $(A)$ and $(C)$ are nonempty. It then follows that $(B')$ can only violate \textbf{rule 6}. When $(A|B|C)$ violates \textbf{rule 4}, we then have
\begin{align*}
    \Pot(A|B|C) &\leq 1+ \tfrac{1}{2}\sign^-(s_{i-k} \cdot s_{i-k-1}),\\
    \Pot(A')&=-\Pot(A)-2\leq -4,\\
    \Pot(C')&=-\Pot(C)-2\leq 1+\tfrac{1}{2}\sign^+(s_{i-k}\cdot s_{i-k-1}),\\
    \Pot(B') &\leq -4,\\
    \Pot(A'|B'|C')&\leq -7+\tfrac{1}{2}\sign^+(s_{i-k}\cdot s_{i-k-1}).
\end{align*}
But $\Pot(A'|B'|C') \geq -5- \tfrac{1}{2}\sign^-(s_{i-k} \cdot s_{i-k-1})$ which is impossible. Similarly, when $(A|B|C)$ violates \textbf{rule 5}, we have
\begin{align*}
    \Pot(A|B|C) &\leq -5- \tfrac{1}{2}\sign^+(s_{i-k} \cdot s_{i-k-1}),\\
    \Pot(A')&=-\Pot(A)-2\leq 3,\\
    \Pot(C')&=-\Pot(C)-2\leq 1+\tfrac{1}{2}\sign^+(s_{i-k}\cdot s_{i-k-1}),\\
    \Pot(B') &\leq -4,\\
    \Pot(A'|B'|C')&\leq \tfrac{1}{2}\sign^+(s_{i-k}\cdot s_{i-k-1}).
\end{align*}
But $\Pot(A'|B'|C') \geq 1+ \tfrac{1}{2}\sign^+(s_{i-k} \cdot s_{i-k-1})$ which is impossible. When $(C)$ is empty, the above argument works identically, that is,
\begin{align*}
    \Pot(A|B) &\leq 1+ \tfrac{1}{2}\sign^-(s_{i-k} \cdot s_{i-k-1}),\\
    \Pot(A')&=-\Pot(A)-2\leq -4,\\
    \Pot(B') &\leq -3-\tfrac{1}{2}\sign^+(s_{i-k}\cdot s_{i-k-1}),\\
    \Pot(A'|B')&\leq -7-\tfrac{1}{2}\sign^+(s_{i-k}\cdot s_{i-k-1}),\\
    \Pot(A'|B') &\geq -5- \tfrac{1}{2}\sign^-(s_{i-k} \cdot s_{i-k-1})
\end{align*}
for \textbf{rule 4}, and for \textbf{rule 5}, we have that
\begin{align*}
    \Pot(A|B) &\leq -5- \tfrac{1}{2}\sign^+(s_{i-k} \cdot s_{i-k-1}),\\
    \Pot(A')&=-\Pot(A)-2\leq 3,\\
    \Pot(B') &\leq -3-\tfrac{1}{2}\sign^+(s_{i-k}\cdot s_{i-k-1}),\\
    \Pot(A'|B')&\leq -\tfrac{1}{2}\sign^+(s_{i-k}\cdot s_{i-k-1}),\\
    \Pot(A'|B')&\geq 1+\tfrac{1}{2}\sign^+(s_{i-k} \cdot s_{i-k-1}).
\end{align*}

When $(A)$ is empty, $(B')$ violates \textbf{rule 4} if $(B|C)$ violates \textbf{rule 5} and vice versa. Thus, when $(B|C)$ violates \textbf{rule 4},
\begin{align*}
    \Pot(B|C) &\leq 1+ \tfrac{1}{2}\sign^-(s_{i-k} \cdot s_{i-k-1}),\\
    \Pot(C')&\geq-3-\tfrac{1}{2}\sign^+(s'_{i-k}\cdot s'_{i-k-1}),\\
    \Pot(C)&\geq 1+\tfrac{1}{2}\sign^+(s'_{i-k}\cdot s'_{i-k-1}),\\
    \Pot(B') &\leq -5,\\
    \Pot(B) &\geq 3,\\
    \Pot(B|C)&\geq 4+\tfrac{1}{2}\sign^+(s'_{i-k}\cdot s'_{i-k-1}).
\end{align*}
Similarly, when $(A|B|C)$ violates \textbf{rule 5}, we have that
\begin{align*}
    \Pot(B|C) &\leq -5- \tfrac{1}{2}\sign^+(s_{i-k} \cdot s_{i-k-1}),\\
    \Pot(C')&\geq-3-\tfrac{1}{2}\sign^+(s'_{i-k}\cdot s'_{i-k-1}),\\
    \Pot(C)&\geq 1+\tfrac{1}{2}\sign^+(s'_{i-k}\cdot s'_{i-k-1}),\\
    \Pot(B') &\leq 1,\\
    \Pot(B) &\geq -3,\\
    \Pot(B|C)&\geq -2+\tfrac{1}{2}\sign^+(s'_{i-k}\cdot s'_{i-k-1})
\end{align*}
which is impossible.
	\begin{figure}[ht]
		\includegraphics[scale=0.5]{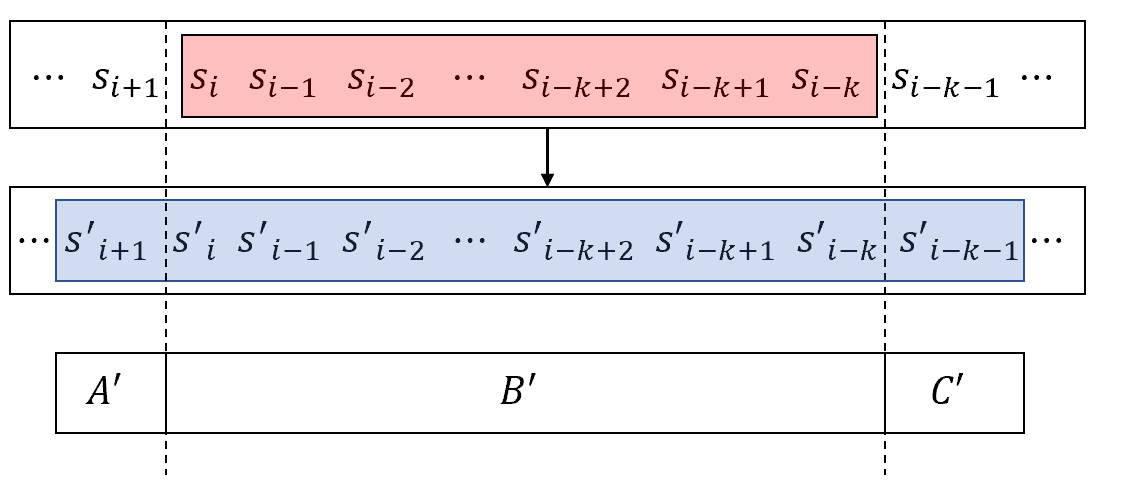}
		\caption{New minimal poison subword cannot contain a subword corresponding to the previous minimal poison subword.}
	\end{figure}
\vspace{0.5cm}
\item The subword $(s'_i,\dotsc,s'_{i-k})$ is contained in the subword $(s'_j,\dotsc,s'_{j-l})$.\par
    Suppose both $(B)$ and $(A'|B'|C')$ violates \textbf{rule 6}. Since $(B)$ is the rightmost, we have that $$\Pot(C)\geq -1-\tfrac{1}{2}\sign^+(C\cdot s'_{j-l}),$$ and therefore, $$\Pot(C')\geq 1-\tfrac{1}{2}\sign^+(C'\cdot s_{j-l}).$$ Thus, $(A'|B'|C')$ violates \textbf{rule 6}, and therefore, we have that $$\Pot(A'|B'|C')\leq -2-\sign^+(s_{j+1}\cdot A)-\tfrac{1}{2}
    \sign^+(C'\cdot s_{j-l}).$$ This gives us a new bound $\Pot(A'|B')\leq-3-\sign^+(s_{j+1}\cdot A)$ which is impossible since $(A|B|C)$ is a minimal poison word. The case $(A'|B'|C')$ when violates the \textbf{rule 4} or \textbf{5} can be worked out similarly.
    
    When $(B)$ violates \textbf{rule 4} or \textbf{rule 5}, then automatically $(A')$ is empty. The potential of $(C)$ is bounded and so is $\Pot(C')$. Combining this with $\Pot(B')$ from $\Pot(B)$, we conclude that $\Pot(B'|C')$ cannot be small enough to violate \textbf{rule 5} or \textbf{rule 4}.
	\vspace{0.5cm}
\item The subword $(s'_j,\dotsc,s'_{j-l})$ is on the right of the subword $(s'_i,\dotsc,s'_{i-k})$.\par
	\begin{figure}[ht]
		\includegraphics[scale=0.5]{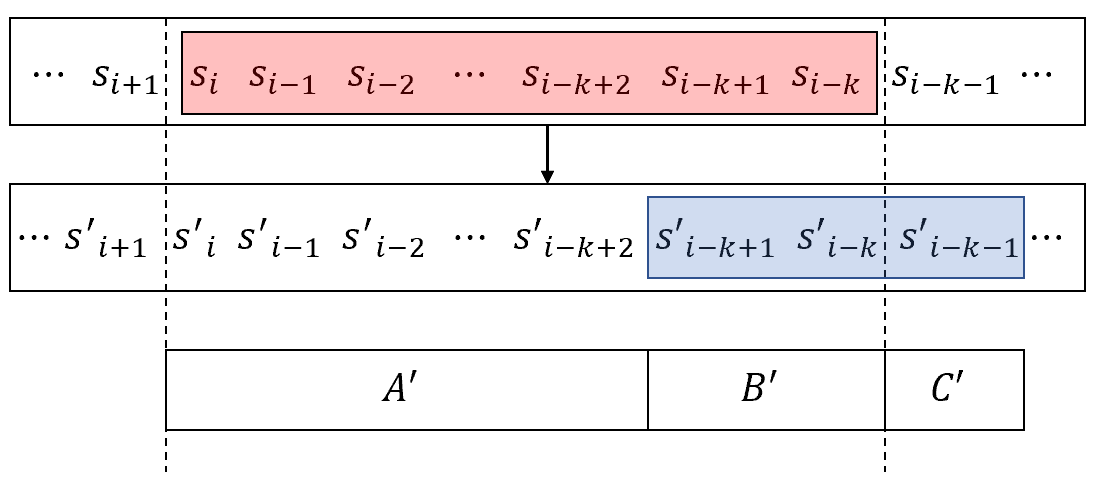}
		\caption{New minimal poison subword cannot appear on the right.}
	\end{figure}
	
Now assume that $(s'_i,\dotsc,s'_{i-k})$ and $(s'_j,\dotsc,s'_{j-l})$ intersect without each containing each other. Suppose $i\geq j\geq i-k\geq j-l$. Let the intersection be $(B')$, $(s'_i,\dotsc,s'_{i-k})=(A'|B')$, and $(s'_j,\dotsc,s'_{j-l})=(B'|C')$. We show $(B')$ is empty and use the fact that $(A|B)$ is the rightmost minimal poison subword to get a contradiction. In this case, the rewriting is given by
\begin{align*}
	&(s_{i+1}|A|B|s_{i-k-1})-\sign(s_i)(-1,2k,2k-1,\dotsc,2k-1,2k,-1)\\
	&=(s'_{i+1}|A'|B'|s'_{i-k-1}).
\end{align*} All coefficients in $(A)$ and $(B)$ have the same sign with all coefficients in $(C)$ the opposite sign. $(A'|B')$ violates \textbf{rule 6}. Thus,
\begin{align*}
	\Pot(A|B)&\leq -2-\sign^+(s_i\cdot s_{i+1}),\\
	\Pot(A)&\geq-2- \sign^+(s_i\cdot s_{i+1}) & \text{\textbf{(rule 6})},\\
	\Pot(B)&\geq -2 & \text{(\textbf{rule 6})},\\
	\Pot(C)&\geq -1-\tfrac{1}{2}\sign^+(s_{i-k-1}\cdot s_{j-l-1}) & \text{(\textbf{rule 6})},\\
	\Pot(B'|C')&\leq -3-\tfrac{1}{2}\sign^+(s'_{j-l}\cdot s_{j-l-1}) & \text{(Violation of \textbf{rule 6})},\\
	\Pot(A')&=-\Pot(A)-2\leq \sign^+(s_i\cdot s_{i+1}) & \text{(Rewriting)},\\
	\Pot(A'|B')&\geq-2+\sign^+(s_i\cdot s_{i+1}) & \text{(Rewriting)}.\\
	\intertext{By combining the inequalities for $\Pot(A')$ and $\Pot(A'|B')$, we get}
    \Pot(B')&\geq -2.
	\intertext{Hence, combining the inequalities for $\Pot(B')$ and $\Pot(B'|C')$ allows us to write}
	\Pot(C')&\leq -1-\tfrac{1}{2}\sign^+(s'_{j-l}\cdot s_{j-l-1}).
\end{align*}
However, it then follows that $\Pot(C)\leq -3-\tfrac{1}{2}\sign^+(s'_{j-l}\cdot s_{j-l-1})$ which contradicts the previous bound on $\Pot(C)$. Thus, $(B)$ is empty. Similarly, we can show that $(B)$ is empty assuming $(A|B)$ violates either \textbf{rule 4} or \textbf{rule 5}.

Since $(B)$ is empty, two consecutive poison subwords are disjoint, and thus, we can directly compare the potentials. $(C')$ is a minimal poison subword. Therefore, we have that $\Pot(C')\leq -2-\sign^+(c_{j}\cdot c_{j+1})-\tfrac{1}{2}\sign^+(c_{j-l}\cdot c_{j-l-1})$. Because $(A)$ and $(C)$ has opposite sign, the rewriting rule suggests that $\Pot(C)=\Pot(C')-2$, which means that $(C)$ already violates \textbf{rule 6} contradicting the fact $(A)$ is the rightmost minimal poison subword. \qedhere
\end{enumerate} \end{proof}
We point out that the proof does not require the minimal poison subword to be rightmost except for the last two cases. The purpose of choosing the rightmost one is so that one can rule out other minimal poison subwords from appearing at certain places. Modifying the proof in the same spirit, we have the following corollary.

\begin{cor}\label{disjoint}
If a word has only one minimal poison subword, then all consecutive rewriting rules are disjoint. In particular, if $P$ is a $n$-reduced polynomial, then all consecutive rewriting rules for $P+1$ are disjoint.
\end{cor}
If the only minimal poison subword comes from \textbf{rule 4} or \textbf{rule 5}, then this follows by the first two cases in the proof of the previous theorem. If it violates \textbf{rule 6}, then this follows from the last case of the proof applied to the left side instead of the right side.

We now show that there exists a $n$-minimal representative of given element $g\in G$ whose polynomial part and principal part are $n$-reduced. We start with the following lemma.

\begin{lemma}\label{first_3_rule_automatic}
Suppose $F(X) = \sum_{i= \ell}^m c_i X^i$ is a $n$-minimal representative of $g=(x,t^n)$. Then polynomial and principal parts of $F(X)$ satisfy \textbf{rules 1, 2} and \textbf{3}.
\end{lemma}
If not, by adding or subtracting $X^d(X^2-(2k+1)X+1)$ with some appropriate choice of $d$, one can obtain a new Laurent polynomial with strictly less $n$-length, contradicting the fact $F(X)$ is a $n$-minimal representative of $g$.

\begin{proposition}\label{rewriting}
Let $F(X) = \sum_{i=\ell}^m c_i X^i$ be a $n$-minimal representative of $g=(x,t^n)$. Let $r$ be an integer, and suppose that $\ell < 0 \leq r \leq m$. Then there exists a deterministic process that takes $F(X)$ to the polynomial $G(X)$ given by
$$
G(X) = \sum_{i=r}^{M} B_i X^i + (B_{r-1} + c_{r-1})X^{r-1} + \sum_{i = \ell}^{r-2} c_i X^i
$$
such that $G(X)$ is a $n$-minimal representative of $x$ and where $\sum_{i=r}^m B_i X^i$ is $n$-reduced. Moreover,  we have that $B_{r-1} \in \{-1, 0, 1\}$, and if $B_{r-1} = \pm 1$, then $B_{r} B_{r-1} < 0$ and $c_{r-1} B_{r}  \geq 0$. We have that $B_r \in \{ c_r, c_r \pm 1, c_r \pm 2k, c_r \pm (2k+1)\}$. Finally, $M = m$ or $M = m \pm 1$.

\begin{proof}
We apply the rewriting procedure as in Theorem \ref{unique.write} on $F(X)$. Start with $c_m X^m$. We can only check local rules, and by the lemma, $c_m X^m$ is $n$-reduced. Thus, $C_m=B_M$ and $m=M$. Similarly, for $F_r(X)=c_m X^m + \cdots + c_{r} X^{r}$, if it is $n$-reduced, we are done and $m=M$. If not, by the lemma, $F_r(X)$ cannot violate local rules, and thus, global rules apply. As in Theorem \ref{unique.write}, we start rewriting from the rightmost minimal poison subword of $(c_m,\dotsc,c_r)$. All rewriting rules after the initial rewriting are disjoint from $(c_r)$ and clearly we have $B_{r-1} \in \{-1, 0, 1\}$ and $B_r \in \{ c_r, c_r \pm 1, c_r \pm 2k, c_r \pm (2k+1)\}$ depending on whether $c_r$ is contained in the rightmost minimal poison subword. Finally, $M=m+1$ or $M=m-1$ depending on whether $c_m$ is contained in the rightmost minimal poison subword in the sequence of rewriting rules. But once $c_m$ is contained in one of the rightmost minimal poison subword, the sequence of rewriting rules terminates.
\end{proof}

\end{proposition}

\begin{theorem}\label{decomposition.thm}
Suppose that $g= (x,t^n)$ is an element of $G$. Then $x$ has an $n$-minimal representative whose polynomial and principal part are both $n$-reduced.
\end{theorem}
We show that we can modify any $n$-minimal representative of $x$ until the desired form is achieved. By Lemma  \ref{first_3_rule_automatic}, we have any $n$-minimal representative must satisfy \textbf{rule 1, 2} and \textbf{3}. Thus, Proposition \ref{rewriting} can be applied on both polynomial and the principal part. To show that this is sufficient, we need the following proposition.

\begin{proposition}\label{truncate}
Let $n \in \Z$.
\begin{enumerate}
\item If $P(X) = \sum_{i=0}^m c_i X^i$ is $n$-reduced, then so is $\sum_{i=j}^m c_i X^k$.
\item Let $P(X)  \in \Z[X]$. Then $P(X)$ is $n$-reduced if and only if $XP(X)$ is $(n+1)$-reduced.
\end{enumerate}
\end{proposition}
\begin{proof}
For the first statement, we let $P_j(X) = \sum_{i=j}^m c_i X^k$. We claim that $P_j$ is $n$-reduced, and towards that end, we may represent $P_j$ as the string $(c_m,\dotsc, c_j,\\ 0,\dotsc, 0)$ where there are $j-1$ zero's at the end of the string. Since $P$ is $n$-reduced and $(c_m,\dotsc, c_j,)$ is a substring of $(c_m,\dotsc, c_0)$, it is easy to see that the 
\textbf{rules 1, 2}, and \textbf{3} are satisfied for coefficients with indices between $j$ and $m$. Since the coefficients of $P_j$ with index less than $j$ are zero, we have that the \textbf{rules 1}, \textbf{2}, and \textbf{3} are all satisfied.

Let $m \geq t \geq i \geq j$, and suppose that $Q = (c_t,\dotsc, c_i)$ is a poison subword with respect to \textbf{rules 4, 5}, or \textbf{6}. If $i>j$, then $(c_t,\dotsc, c_i)$ is a poison subword of $P$ with respect to \textbf{rules 4, 5}, or \textbf{6} which is a contradiction. Thus, if $P_j$ contains a poison subword, it must contain the coefficient $c_j$.

Suppose that $(c_m,\dotsc, c_i)$ violates \textbf{rule 4}. We then would have that $$(c_m,\dotsc, c_j) = \pm(1, c_{m-1}^\prime,\dotsc, c^\prime_j)$$ where $c_{m-i^\prime}^\prime \in \{-k+1,-k\}$ for $i^\prime \neq i$ and $c_{m-i}^\prime \in \{-k,-k-1\}$. We have that $\sign^\prime(c_m \cdot c_{j-1}) = \sign^\prime(c_m \cdot 0) = \sign^\prime(0)$.  We then have $\text{Pot}(c_m,\dotsc,c_i) \leq 0$. However, since $(c_m,\dotsc, c_i)$ is a substring of $P$, we have that $(c_m,\dotsc,c_i)$ satisfies \textbf{rule 4}. In particular, if $\sign(c_m \cdot c_{i}) > 0$, then $\text{Pot}(c_m,\dotsc, c_{i}) \geq 2$ which is a contradiction. If $\sign(c_m \cdot c_{i-1}) \leq 0$, then $\text{Pot}(c_m,\dotsc, c_i) > 1$ which is also a contradiction. Thus, $P_j$ satisfies \textbf{rule 4}.

Suppose that $(c_m,\dotsc, c_i)$ violates \textbf{rule 5}. We would have that $c_t \in \{k+1, k+2\}$ where $c_s \in \pm \{k-1,k\}$ for $s \neq i$ and $c_{i} \in \{k, k+1\}$. Since $c_{i-1} = 0$, we have that $\sign(c_m \cdot c_{i-1}) = 0$. Thus, we have that $\text{Pot}(c_m,\dotsc, c_{i-1}) <-6$. Since $(c_m,\dotsc, c_i)$ is a substring of $(c_m,\dotsc, c_0)$, we have that $(c_m,\dotsc, c_i)$ satisfies \textbf{rule 5}. Hence, if $\sign(c_m \cdot c_{i-1}) \geq 0$, then $\text{Pot}(c_m,\dotsc, c_i) > -6$ which is a contradiction. If $\sign(c_m \cdot c_{i-1}) < 0$, then $\text{Pot}(c_t,\dotsc, c_i) > -4$ which is also a contradiction. Thus, $P_j$ must satisfy rule \textbf{rule 5}.

Finally, assume that $(c_t,\dotsc, c_i)$ violates \textbf{rule 6}. We then have that $c_s \in \{k,k+1\}$ for $s \in \{i,j\}$ and $c_s \in \pm\{k-1,k\}$ otherwise. By definition of \textbf{rule 6}, we have that $\text{Pot}(c_t,\dotsc, c_{i-1}) < -3 - \sign^+(c_t \cdot c_{t+1})$. However, $(c_t,\dotsc, c_i)$ is a substring of $P$ which implies that it satisfies \textbf{rule 6}. If $\sign(c_t \cdot_{i-1}) \geq 0$, then $\text{Pot}(c_j, \cdot, c_i) > - 3 - \sign^+(c_t \cdot c_{t+1})$ which is a contradiction. If $\sign(c_t \cdot c_{t+1}) < 0$, then $\text{Pot}(c_t,\dotsc, c_{j}) \geq -1 - \sign^+(c_t \cdot c_{t+1}) > 0$ which is also a contradiction. Thus, $P_j$ must satisfy \textbf{rule 6}.

For the second statement, we may proceed using similar arguments as for the first statement.
\end{proof}

\begin{proof}[Proof of Theorem \ref{decomposition.thm}]
    Suppose that $F(X)$ represents $g=(x,t^n)$ with both non-trivial polynomial and principal parts. Applying Proposition \ref{rewriting} on the principal part $F(X)$ first, there exists $F'(X)=\sum_{i=l<0}^{M} B_i X^i$ that represents $g$. Now we apply Proposition \ref{rewriting} on the polynomial part and we have $$G(X) = \sum_{i=0}^{M} C_i X^i + (C_{-1} + B_{-1})X^{-1} + \sum_{i = \ell}^{-2} B_i X^i.$$ If $C_{-1}=0$, then we are done. If not, we have $C_0 C_{-1}<0$ and $C_0 B_{-1}\geq 0$. Observe that by Proposition \ref{truncate}, $\sum_{i = \ell}^{-2} B_i X^i$ is $n$-reduced. By Corollary \ref{disjoint}, we have that all rewriting rules for the principal parts are disjoint from the polynomial part except $C_0$. Rewriting using Proposition \ref{rewriting} on the principal part, we obtain a Laurent polynomial that is $n$-reduced on the principal part, and $C_0 C_{-1}<0$ guarantees that the polynomial part is $n$-reduced as well.
\end{proof}

\section{Stability of $n$-reduced polynomials}

In this section, we attempt to count the number of all $n$-reduced polynomials. Intuitively, we start with $P=0$ which is trivially $n$-reduced and keep adding 1 repeatedly until we fail to have a $n$-reduced polynomial. By quantifying this failure, we classify all $n$-reduced polynomials with the nonnegative leading coefficient. We start with a basic lemma.
\begin{lemma}
Suppose that $P=(\cdots,c_1,c_0)$ is a $n$-reduced polynomial. Then $P+1$ is $n$-reduced or the rightmost minimal poison subword of $P+1$ contains either $c_1$ or $c_0+1$.
\end{lemma}
If there is no minimal poison subword of $P+1$ containing any of the coefficients, then $P+1$ is $n$-reduced.

\begin{definition}Suppose that $P=(\cdots,c_1,c_0)$ is a $n$-reduced polynomial where $P+1$ is not $n$-reduced. If the rightmost minimal poison subword contains $c_1$ but not $c_0$, we say $P+1$ fails to be $n$-reduced by a \textbf{sign change violation}. If the rightmost minimal poison subword contains $c_0$, we say $P+1$ fails to be $n$-reduced by a \textbf{potential value change violation}.
\end{definition}

We see that in order for a polynomial $P+1$ to have a sign change violation, we must have $c_0=0$ and that there exists an upper bound on the potential of the rightmost minimal poison subword. We then have that the polynomial $P+1$ is no longer $n$-reduced due to the change in sign of $c_0$ on the condition of the rules. Similarly, suppose $P+1$ fails to be $n$-reduced by a potential change violation and $A'$ be its rightmost minimal poison subword. Let $A$ be the subword on the exact same location before adding 1. The potential of $A'$ and the potential of $A$ differ by $-2$. Since $A'$ violates one of the rules, the potential of $A$ is bounded above. In both cases, we have an upper and lower bound on the potential with different possible rightmost minimal poison subwords. Thus, we are able to classify all cases.

\renewcommand{\labelenumi}{\arabic{enumi})}
\renewcommand{\labelenumii}{(\roman{enumii})}
\begin{proposition}\label{classification}
Suppose that $P$ is a $n$-reduced polynomial of degree $m$ and $P+1$ is not $n$-reduced. Then, $P$ falls into one of these categories:
\begin{itemize}[topsep=4pt, itemsep=4pt, leftmargin=*]
\item {Sign change violations}

\begin{enumerate}[topsep=4pt, itemsep=4pt]
\item $P+1$ violates \textbf{rule 1}.\\In this case, $P=(-k-2,0)$ and $n\leq 1$.
\item $P+1$ violates \textbf{rule 2}.\begin{enumerate}
    \item $P=(-k-1,0)$ and $n>1$.
    \item $P=(\cdots,c_2,-k-1,0)$ and $c_2<0$.
    \item $P=(\cdots,c_2,-k,0)$ and $c_2\geq 0$.
\end{enumerate}
\item $P+1$ violates \textbf{rule 3}.\\In this case, $P=(1,-k+1,0)$ and $n\leq 1$.

\item $P+1$ violates \textbf{rule 4}.\\For $P=(c_{m-1},\dotsc,c_1,0)$, $m>n$, $c_1=-k$ and $$\Pot(1,c_{m-1},\dotsc,c_1)=1.$$

\item $P+1$ violates \textbf{rule 5}.\\For $P = (c_m,\dotsc, c_1, 0)$, $m\geq n$, $c_1=-k$ and $$\Pot(c_m,\dotsc, c_1)=-5.$$

\item If $P+1$ violates \textbf{rule 6}.\\Let $P = (\cdots |A|0)$ and $(A)=(c_j,\dotsc, c_1)$ be the minimal poison subword of $P+1$. Then, $c_1=-k$ and $$\Pot(A)=-2-\sign^+(c_j \cdot c_{j+1}).$$
\end{enumerate}

\item Potential value change violations

\begin{enumerate}[topsep=4pt, itemsep=4pt]
\item $P+1$ violates \textbf{rule 1}.\\In this case, $P=k+2$ and $n\leq 0$.
\item $P+1$ violates \textbf{rule 2}.
    \begin{enumerate}
        \item $P=k+1$, $n>0.$
        \item $P=(\cdots,c_1,k+1)$, $c_1\geq 0.$
        \item $P=(\cdots,c_1,k)$, $c_1<0.$
    \end{enumerate}
\item $P+1$ violates \textbf{rule 3}.\\In this case, $n\leq 0$, $P=(-1,k-1)$.
\item $P+1$ violates \textbf{rule 4}.\\For $P=(-1,c_{m-1},\dotsc, c_0)$, $m>n$
\begin{enumerate}
    \item $c_0=k$ and $$\Pot(c_{m-1},\dotsc,c_1)\in \{2,3\}.$$
    \item $c_0=k-1$ and $$\Pot(c_{m-1},\dotsc,c_1)=1.$$
\end{enumerate}

\item If $P+1$ violates \textbf{rule 5}.\\For $P = (c_m,\dotsc, c_1, c_0)$, $m\geq n$
\begin{enumerate}
\item $c_0=k-1$ and $$\Pot(c_m,\dotsc, c_1) = -5.$$
\item $c_0=k$ and $$\Pot(c_m,\dotsc, c_1) \in \{-4, -3\}.$$
\end{enumerate}
\item If $P+1$ violates \textbf{rule 6}.\\Let $(A|c_0)=(c_j,\dotsc, c_1, c_0)$ be the minimal poison subword of $P+1$. Then,
\begin{enumerate}
\item $c_0=k-1$ and $$\Pot(A) = -2 - \sign^+(c_j \cdot c_{j+1}).$$
\item $c_0=k$ and $$\Pot(A) \in \{ -1 - \sign^+(c_j \cdot c_{j+1}),-\sign^+(c_j \cdot c_{j+1})\}.$$
\end{enumerate}
\end{enumerate}

\end{itemize}
\end{proposition}
We note that some symmetries are present in these classifications. For any $P=QX+1$ that violates any of the rules for given $n$ by sign change, $-Q+1$ violates the same rule for $n-1$ by potential value change, which reduces the number of cases to consider. Furthermore, since all rewriting rules for $P+1$ are disjoint when it is not $n$-reduced, this suggests that we can break $P+1$ into smaller pieces with specific potential condition. Moreover, since the associated rewriting starts with $(1,-k,-k+1,\cdots,-k)$, each piece either falls into one of the above categories. That means we can inductively classify all of $n$-reduced polynomials what we will later call $n$-\textbf{types} and $n$-\textbf{classes}.

\section{$n$-Types and $n$-Classes}
%\begin{definition}
%$\:$ \\ \noindent Let $P$ be a $n$-reduced polynomial $P$ of degree $m > 0$ written as $(c_m,\dotsc, c_0)$. \par\noindent $P$ is of
%\begin{itemize}[leftmargin=1.5cm]
%\item[\textbf{class $A$}] if \begin{itemize} \item there is some $j < \min(n,m)$ such that $c_i \in \{k+1,k,k-1\}$ for $0 \leq i \leq j$ and $\Pot(c_j,\dotsc,c_0) = -2-\sign^+(c_jc_{j+1})$\end{itemize}.
%\item[\textbf{class $B$}] if
%\begin{itemize}
%\item no subword of $P$ is of type $A$
%\item there exists a $j < \min(n,m)$ such that $c_i \in \{k+1, k, k-1\}$ for $0 \leq i \leq j$ and $\Pot(c_j, 
%\cdots, c_0) = - 1 - \sign^+(c_j \cdot c_{j+1})$.
%\end{itemize} 
%\item[\textbf{class $C$}] if
%\begin{itemize}
%\item no subword of $P$ is of type $A$ or $B$
%\item there exists a $j < \min(n,m)$ such that $c_i \in \{k+1, k, k-1\}$ for $0 \leq i \leq j$ and $\Pot(c_j,\dotsc, c_0) = -\sign^+(c_jc_{j+1})$.
%\end{itemize}
%\item[\textbf{class $E$}] if
%\begin{itemize}
%\item $P = (1,c_{m-1},\dotsc,c_0)$ such that $c_0 =-1$ and %$\Pot(c_{m-1},\dotsc,c_1)=1$.
%\end{itemize}
%\end{itemize}
%\end{definition}

In this section, we give the precise definition of $n$-types and $n$-classes based on the classification in Proposition \ref{classification}. The type of $P$ captures the behavior of $P+1$, and the class of $P$ captures the behavior of $XP+C$. As a result, the definition of the class of $P$ will heavily depend on the potential of a final segment of $P$. Before giving the formal definitions, we give the intuition behind them.

As mentioned in Proposition \ref{classification}, there is a symmetry between sign change violations and potential value change violations. Hence, there's a reduction in cases. For later use, we classify polynomials whose leading coefficient is nonnegative. Thus, the sign change violations by rule 4 will be preferred over potential value change violations by rule 4.

We start with the generic case. If both $P$ and $P+1$ are $n$-reduced and the constant coefficient of $P$ is non-negative, we say $P$ is of $n$-interior type. If the constant coefficient of $P$ is negative, then it follows that $P+1$ is also $n$-reduced, and we say $P$ is of $n$-negative type. All other polynomials, i.e.\ when $P+1$ is not $n$-reduced, are said to be of $n$-boundary type. Geometrically, this is dividing the ball into three different regions - the interior, boundary and negative with respect to the direction $(+1)$.

Now we discuss $n$-classes. We begin with the \textbf{stable} $n$-class, shortly, $n$-class $S$. This class consists of $n$-reduced polynomials $P(X)$ such that for any positive $C$, the rightmost minimal poison subword of $X P(X)+C$ is disjoint from the subword corresponding to $X P(X)$. Necessarily, $P+1$ is $n$-reduced whenever $P$ is of $n$-class $S$. These are the most common $n$-reduced polynomials. The class splits into 3 separate sub-classes $S_{-}$, $S_{0}$ and $S_{+}$ depending on the sign of the constant coefficient.

\begin{table}[ht!]
    \centering
    \begin{tabular}{c|c|c}
        \textit{Violations} & Corresponding $n$-class & Related $n$-classes \\\hline
        Sign Change, \textbf{rule 4} & $E_1$ & $E_2$, $E_3$ \\\hline
        Potential Change, \textbf{rule 5} & $U_{-3}^t$, $U_{-4}^t$, $U_{-5}^t$ & \\\hline
        Potential Change, \textbf{rule 6} & $U_{0}$, $U_{-1}$, $U_{-2}$ & \\\hline
    \end{tabular}
    \vspace{0.2cm}\\ Generic classes: $S_{-}$, $S_{0}$, $S_{+}$
\end{table}

The opposite of the stable $n$-class is the \textbf{unstable} $n$-class which we write as $n$-class $U$. These are the polynomials $P(X)$ where the rightmost minimal poision word of $XP(X)+C$ violates some global rule for some positive $C$. These violations are classified in Proposition \ref{classification}.

We first consider $n$-reduced polynomials $P$ such that $XP(X)+k$ or $XP(X)+k+1$ is not $(n+1)$-reduced and violates \textbf{rule 6} as we saw in Proposition \ref{classification}. This class splits into sub-classes $U_{0}$, $U_{-1}$ and $U_{-2}$ by the truncated potential of the rightmost minimal poison subword $- \sign^+(c_j \cdot c_{j+1}),$ $-1- \sign^+(c_j \cdot c_{j+1}),$ and $-2 - \sign^+(c_j \cdot c_{j+1})$, respectively. We note that $P'(X)= XP(X)+k-1$ or that $P'(X) = XP(X)+k$ if $(n+1)$-reduced is of $(n+1)$-class $U$ by the same logic. Furthermore, if we do the same analysis to $P'(X)$, we can control what is the truncated potential of the rightmost minimal poision subword of $XP'(X)+C'$. For example, if $P$ is of class $U_{0}$, then $XP+k+1$ is not $(n+1)$-reduced and $XP+k$ will be class $U_{-1}$ as the potential of $k$ is $-1$. Also note that since the sign change and potential violations of rule 6 share the same ``edge cases'', the three sub-classes will suffice.

Similarly, we will define $n$-classes $U^t_{-2}$, $U^t_{-3}$ and $U^t_{-4}$ for \textbf{rule 5} in Proposition \ref{classification} where the subscript denotes the potential. Finally, for \textbf{rule 4}, we will define $n$-classes $E_1$, $E_2$, and $E_3$. Note that the $U^t$ classes and the $E$ classes are dual to each other and are both associated to top rules.

Based on these intuition and ideas, we give the following formal recursive definition:
\begin{table}[ht!]
\centering
  \begin{tabular}{ c | c | c  c | c }
  \hline
  Type & $n$-class of $Q$ & \multicolumn{2}{c|}{$Q=dX+c$} &$n$\\ \hline\hline
  $n$-initial type & - & \multicolumn{2}{c|}{$0$}& $n\leq 0$\\ \hline\hline
  \multirow{5}{*}{$n$-interior type} & {$S_{0}$} & $d=0$&$c=0$ & $n> 0$\\\cline{2-5}
  & $S_{+}$ & $d=0$ & $0<c\leq k$ & $n\le 0$\\\cline{2-5}
  & $U^t_{-3}$ & $d=0$ & $c=k+1$ & $n\leq 0$\\\cline{2-5}
  & $S_{+}$ & $d=0$&$0<c\leq k-1$ & $n> 0$\\\cline{2-5}
  & $U_{0}$ & $d=0$&$c=k$ & $n> 0$\\ \hline
  \multirow{2}{*}{$n$-negative type} & $E_1$ & \multicolumn{2}{c|}{$X-k+1$} & $n\leq 0$ \\\cline{2-5}
   & $E_3$ & \multicolumn{2}{c|}{$X-k+2$} & $n\leq 0$ \\ \hline
  \multirow{2}{*}{$n$-boundary type} & $U^t_{-5}$ & $d=0$&$ c=k+2$ & $n\leq 0$\\\cline{2-5}
    & $U_{-2}$ & $d=0$&$c=k+1$ & $n> 0$\\ \hline
%    & $U_{-2}$ & \multicolumn{2}{c|}{$-X+k+1$} & $n> 0$\\ \hline
%    & $U_{-1}$ & \multicolumn{2}{c|}{$-(k+1)X+k-1$} & $n> 1$\\ \cline{2-5}
%    & $U_{-1}$ & \multicolumn{2}{c|}{$-X^2+(k-1)X+k-1$} & $n> 1$\\ \hline
  \end{tabular}
\end{table}

\begin{definition}[Base case]\label{type_base}
$\:$ \par \noindent Let $Q$ be an integer polynomial with a nonnegative leading coefficient. The $n$-class and $n$-type of $Q$ are given as follows:
\begin{itemize}[leftmargin=0cm]
\item (\textbf{$n$-initial type}) $n \leq 0$, $Q=0$.
\item(\textbf{$n$-interior type}) if $Q=dX+c$, and
    \begin{itemize}[leftmargin=0.5cm]
    \item $n > 0$. Then $Q=0$ has \textbf{$n$-class $S_{0}$}. 
    \item $d = 0,$ $0 < c \leq k$, and $n \leq 0$. Then $Q$ has \textbf{$n$-class $S_{+}$}.
    \item $d=0$, $c = k + 1$, and $n \leq 0$. Then $Q$ has \textbf{$n$-class $U^t_{-3}$}.
    \item $d=0$, $0 < c \leq k -1$, and $n > 0$. Then $Q$ has \textbf{$n$-class $S_{+}$}.
    \item $d=0$, $c = k$, and $n > 0$. Then $Q$ has \textbf{$n$-class  $U_{0}$}.
    \end{itemize}
\item(\textbf{$n$-negative type}) if $Q=dX+c$, and
    \begin{itemize}[leftmargin=0.5cm]
    \item $d=1$ and $c = -k + 1$ is $n$-class $E_1$ if $n \leq 0$.
    \item $d=1$ and $c= -k +2$ is a $n$-class $E_3$ if $n \leq 0$.
    \end{itemize}
\item(\textbf{$n$-boundary type})
    \begin{itemize}[leftmargin=0.5cm]
    \item If $n \leq 0$. Then $Q=k+2$ has \textbf{$n$-class $U^t_{-5}$}.
    \item If $n > 0$. Then $Q=k+1$ has \textbf{$n$-class $U_{-2}$}.
    \end{itemize}
\end{itemize}
\end{definition}

We define types and classes for polynomial with a nonnegative leading coefficient.

\begin{definition}\label{type_induct} $\:$ \\
Suppose that $Q(X) = XP(X) + c$ is a $n$-reduced polynomial whose leading coefficient is positive where $P(X) \neq 0$ and $c$ is constant.

\begin{itemize}[leftmargin=0cm]
\item(\textbf{$n$-interior type}) if 
    \begin{itemize}[leftmargin=0.5cm]
    \item For any $P$ and $0< c \leq k - 2$, $Q$ is of $n$-class $S_{+}$.
    \item $P$ is of any $(n-1)$-class except $E_1$ and $-P$ is not of $(n-1)$-class $U_{-2}$ and $c=0$. Then $Q$ is of $n$-class $S_0$.
    \item $P$ is of $(n-1)$-class $S_{0}$, $S_{+}$, $U_0$, or $U^t_{-3}$ and $c =k-1$. Then $Q$ is of $n$-class $S_{+}$.
    \item $P$ is of $(n-1)$-class $U_{-1}$ and $c = k-1$. Then $Q$ is of $U_{0}$.
%    \item $P$ is of $(n-1)$-class $U_{-3}^t$ and $c= k-1$, then $Q$ is of $n$-class $S_{+}$
    \item $P$ is of $(n-1)$-class $S_{0}$ or $S_+$ and $c =k$. Then $Q$ is of $n$-class $U_{0}$.
    \item $P$ is of $(n-1)$-class $S_{-}$, $E_1$, $E_2$, or $E_3$ and $c= k-1$. Then $Q$ is of $n$-class $U_{0}$.
    \item $P$ is of $(n-1)$-class $U^t_{-4}$ and $c = k-1$. Then $Q$ is of $n$-class $U^t_{-3}$.
\end{itemize}
\item(\textbf{$n$-negative type}) if
\begin{itemize}[leftmargin=0.5cm]
    \item $P$ is $(n-1)$-class $E_1$ and $c = -k+1$. Then $Q$ has $n$-class $E_2$.
    \item $P$ is $(n-1)$-class $E_2$ and $c = -k+1$. Then $Q$ has $n$-class $E_3$.
    \item $P$ is $(n-1)$-class $E_3$ and $c = -k+1$. Then $Q$ has $n$-class $S_{-}$.
    \item $P$ is $(n-1)$-class $E_2$ and $c = -k$. Then $Q$ has $n$-class $E_1$.
    \item $P$ is $(n-1)$-class $E_3$ and $c = -k$. Then $Q$ has $n$-class $E_2$.
%    \item $P$ is $(n-1)$-class $S_{-}$ or $S_{0}$ and $-k-1\leq c < 0$. Then $Q$ is $n$-class $S_{-}$.
    \item $c < 0$. Then $Q$ is $n$-class $S_{-}$ unless otherwise stated above.
\end{itemize}
\item(\textbf{$n$-boundary type (P)}) if
\begin{itemize}[leftmargin=0.5cm]
    \item $P$ is of $(n-1)$-class $U_{-2}$ and $c = k-1$. Then $Q$ is of $n$-class $U_{-1}$.
    \item $P$ is of $(n-1)$-class $U^t_{-5}$ and $c = k-1$. Then $Q$ is of $n$-class $U_{-4}^t$.
    \item $P$ is of $(n-1)$-class $U_{-1}$ and $c = k$. Then $Q$ is of $n$-class $U_{-2}$.
    \item $P$ is of $(n-1)$-class $U^t_{-4}$ and $c = k$. Then $Q$ is of $n$-class $U^t_{-5}$.
    \item $P$ is of $(n-1)$-class $U_{0}$ and $c = k$. Then $Q$ is of $n$-class $U_{-1}$.
    \item $P$ is of $(n-1)$-class $U^t_{-3}$ and $c = k$. Then $Q$ is of $n$-class $U^t_{-4}$.
    \item $P$ is of $(n-1)$-class $S_{-}$, $E_1$, $E_2$, or $E_3$ and $c = k$. Then $Q$ is of $n$-class $U_{-2}$.
    \item $P$ is of $(n-1)$-class $S_{+}$ or $S_{0}$ and $c = k+1$. Then $Q$ is of $n$-class $U_{-2}$.
\end{itemize}
\item(\textbf{$n$-boundary type (S)}) if
    \begin{itemize}[leftmargin=0.5cm]
        \item $P$ is of $(n-1)$-class $S_{-}$ such that $-P$ is of $(n-1)$-class $U_{-2}$ and $c = 0$. Then $Q$ is of $n$-class $S_{0}$.
        \item $P$ is of $(n-1)$-class $E_1$ and $c = 0$. Then $Q$ is of $n$-class $S_0$.
    \end{itemize}
\end{itemize}
\end{definition}

\begin{longtable}{ m{2.6cm} | m{3.2cm} | m{3.1cm} | m{2.3cm} }
    \hline
    Type of $Q$ & $n$-class of $XP+c$ & $(n-1)$-class of $P$ & $c$ \\ \hline \hline
    \multirow{7}{*}{$n$-interior} & $S_{+}$ & Any & $0<c\leq k-2$ \\\cline{2-4}
    & $S_{0}$ & except $E_1$, $-U_{-2}$ & $c=0$ \\\cline{2-4}
    & $S_{+}$ & $S_0$, $S_{+}$, $U_0$, $U_{-3}^t$ & $c=k-1$ \\\cline{2-4}
    & $U_{0}$ & $U_{-1}$ & $c=k-1$ \\\cline{2-4}
    & $U_{0}$ & $S_0$, $S_+$ & $c=k$ \\\cline{2-4}
    & $U_{0}$ & $S_-$, $E$ & $c=k-1$ \\\cline{2-4}
	& $U_{-3}^t$ & $U_{-4}^t$ & $c=k-1$ \\
    \hline
    \multirow{4}{*}{$n$-negative} & $E_2$ & $E_1$ & $c=-k+1$ \\\cline{2-4}
    &$E_3$ & $E_2$ & $c=-k+1$ \\\cline{2-4}
    &$S_{-}$ & $E_3$ & $c=-k+1$ \\\cline{2-4}
    &$E_1$ & $E_2$ & $c=-k$ \\\cline{2-4}
    &$E_2$ & $E_3$ & $c=-k$ \\\cline{2-4}
%    &$S_{-}$ & $S_{-}$, $S_{0}$ & $-k-1 \leq c<0$ \\\cline{2-4}
    &$S_{-}$ & except above & $ c<0$ \\
    \hline
    \multirow{8}{*}{$n$-boundary (P)} & $U_{-1}$ & $U_{-2}$ & $c=k-1$ \\\cline{2-4}
    & $U_{-4}^t$ & $U_{-5}^t$ & $c=k-1$ \\\cline{2-4}
    & $U_{-2}$ & $U_{-1}$ & $c=k$ \\\cline{2-4}
    & $U_{-5}^t$ & $U_{-4}^t$ & $c=k$ \\\cline{2-4}
    & $U_{-1}$ & $U_{0}$ & $c=k$ \\\cline{2-4}
    & $U_{-4}^t$ & $U_{-3}^t$ & $c=k$ \\\cline{2-4}
    & $U_{-2}$ & $S_{-}$, $E$ & $c=k$ \\\cline{2-4}
    & $U_{-2}$ & $S_{+}$, $S_{0}$ & $c=k+1$ \\
    \hline
    \multirow{2}{*}{$n$-boundary (S)} & $S_{0}$ & $-U_{-2}$ & $c=0$ \\\cline{2-4}
    & $S_{0}$ & $E_1$ & $c=0$ \\ \hline
\end{longtable}

Although we gave our definitions recursively, it is straightforward to prove inductively that these do correspond to the cases in Proposition \ref{classification}. For example, if $P$ is of $n$-boundary type (P) and $n$-type $U_{-1}$, then:
\begin{itemize}
\item $P = (\cdots|A)$ with $A = (c_j,\ldots,c_0)$ and $\Pot(A) = -1-\sign^+(c_j\cdot c_{j+1})$,
\item $P+1$ is not $n$-reduced and its minimal poison subword is $A' = (c_j,\ldots,c_1,c_0+1)$,
\item $XP+k+1$ will violate \textbf{rule 6} by a potential violation as described in Proposition \ref{classification}.
\end{itemize}

From Proposition \ref{classification}, it is also easy to check that the $n$-classes form a partition of all $n$-reduced polynomials.

\section{Successor function}

In this section, we define the successor of a polynomial to generate all $n$-reduced polynomials with a nonnegative leading coefficient. We will add 1 consecutively to the $n$-reduced polynomial until it is no longer $n$-reduced. When $P(X)+1\in \mathbb{Z}[X]$ is not $n$-reduced, by Proposition \ref{rewriting}, the rewriting rules will give us a new polynomial
\begin{align*}
    Q(X)=\sum_{i=0}^d c_i X^i+c_{-1}X^{-1}
\end{align*}
where $c_{-1}=0$ or $1$. By Proposition \ref{truncate}, we know that $Q(X) -X^{-1} = \sum_{i=0}^d c_i X^i \in \mathbb{Z}[X]$ is also $n$-reduced. Intuitively, we want to define this to be the \textbf{successor} of $P(X)$. However, in some cases, $Q(X)-X^{-1}-1$ is also $n$-reduced, and to make sure that the successor function maps onto the set of $n$-reduced polynomials, we instead define $Q(X)-X^{-1}-1$ to be the successor of $P(X)$ (see (3) in the formal definition).

\begin{definition}\label{def_succ}
Let $P$ be a $n$-reduced polynomial. The \textbf{successor} of $P$, denoted $\S(P)$, is given by the following:
\begin{enumerate}
\item If $P$ is of $n$-initial, $n$-interior, or $n$-negative type, then $\S(P) = P+1$.
\item If $P$ is one of boundary types in the base case, we define the successor as follows:
	\begin{table}[H]\begin{tabular}{c|c}
	 $P$ & $\S(P)$\\ \hline
	 $k+2$ and $n\le 0$ & $X-k+1$ \\
	 $k+1$ and $n> 0$ & $X-k$ \\ \hline
%	 $-X+k+1$ ($n>0$) & $-k-1$ \\
%	 $-X+k+1$ ($n\leq0$) & $-k-2$ \\ 
%	 $-(k+1)X+k-1$ & $X^2-(k-1)X-(k-1)$ \\ 
%	 $-X^2+(k-1)X+k-1$ & $(k-1)X-(k-1)$ \\ \hline
	\end{tabular}\end{table}
\item Suppose that $P$ is of $n$-boundary type (P) and $n$-class $U_{-2}$ or $U_{-5}^t$.. Letting $Q$ be the rewriting of $P+1$ given by Proposition \ref{rewriting}, we then set $\S(P) = Q-1-X^{-1}$.
\item Suppose that $P$ is of $n$-boundary type (P) and $n$-class $U_{-1}$ or $U_{-4}^t$. Letting $Q$ be the rewriting of $P+1$ given by Proposition \ref{rewriting}, we then set $\S(P) = Q-X^{-1}$.
%\item If $P$ is of $n$-boundary type (P) and $n$-class $U_{0}$ or $U_{-3}^t$, then $\S(P) = P+1$.
\item If $P$ is of $n$-boundary type (S), then $\S(P)$ is defined to be the rewriting of $P+1$ given by Proposition \ref{rewriting}.
\end{enumerate} 
If $\S(P) = P+1$, it is called a \textbf{regular successor}. Otherwise, it is called an \textbf{irregular successor}.

\end{definition}

It is clear in the definition that $\S(P)$ is $n$-reduced when it is a regular successor. As an example, we show that $\S(P)$ is $n$-reduced when $P$ is of $U_{-2}$. If $P=XQ(X)+k$ and $Q$ is of $S_{-}$, or $P=XQ(X)+k+1$ and $Q$ is of $S_0$ or $S_+$, then $P+1$ violates \textbf{rule 2}. Thus, it can be checked easily that $\S(P)$ is $n$-reduced. If $P=XQ(X)+k$ and $Q$ is of $U_{-1}$, there exists $d$ such that
$$P=X^{d+1} B+X A + k$$
where $B$ is of stable $n$-type and $(A|k+1)$ is the rightmost minimal poison subword of $P+1$. After rewriting $P+1$ and subtracting $1+X^{-1}$, we have $\S(P)=X^{d+1}(B+1)+XA'-k.$
\begin{figure}[ht!]
	\includegraphics[scale=0.15]{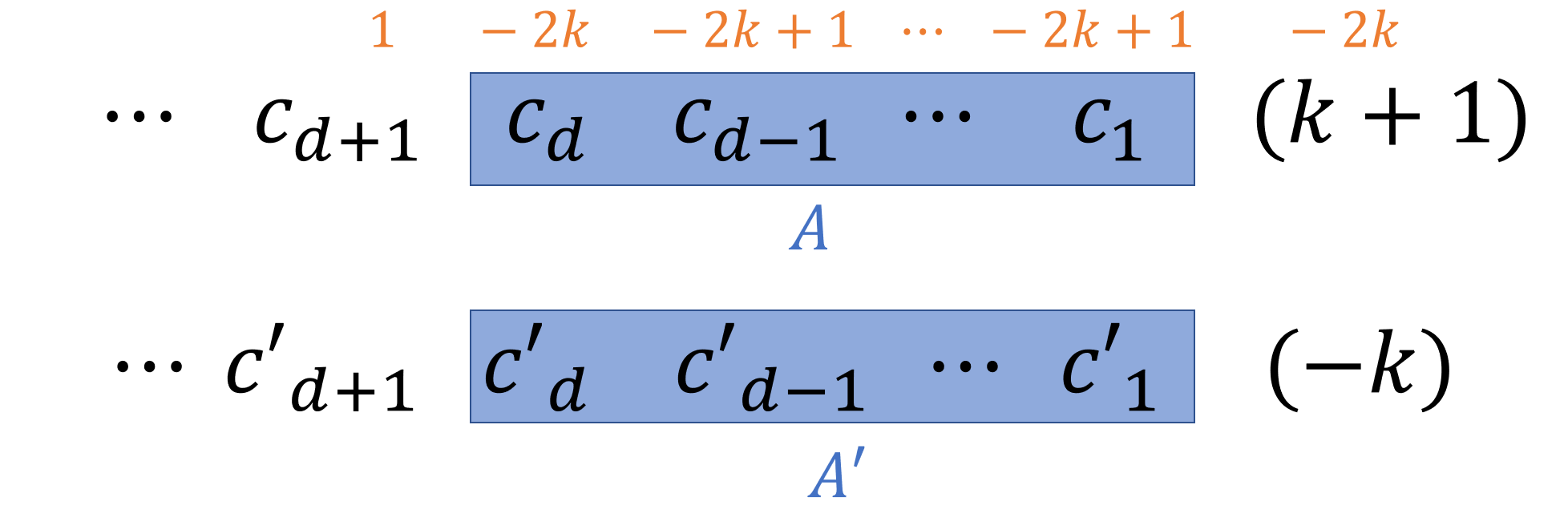}
\end{figure}

Given that
\begin{align*}
	\Pot(A)&=-2-\sign^+(c_d\cdot c_{d+1})\\
	\Pot(A')&=-(-2-\sign^+(c_d\cdot c_{d+1}))-2\\
              &= \sign^+(c_d\cdot c_{d+1}) = -\sign^+(c'_d\cdot c'_{d+1})\\
\end{align*}
and for any subword adjacent to $(-k)$ not containing $c'_d$ we also have $$\Pot(c_d \cdots c_j)>-3-\sign^+(c_d\cdot c_{d+1}),$$ it then follows that 
$$\Pot(c_j \cdots c_1)\leq 0.$$ Thus, $\Pot(c'_j \cdots c'_1)\geq 0$. Since the potential for $A'$ and of any subword of $A'$ is bounded below, we conclude that $\S(P)$ cannot violate \textbf{rule 6}.

For later purposes, we also define the \emph{generalized} successor function $\tilde{\S}$ so that the successor function can be defined for \emph{any} polynomial. This generalized version will be used in this section only. We allow two base cases in the Definition \ref{type_base}:
\begin{itemize}[leftmargin=0cm]
\item(\textbf{$n$-negative type}) if $Q=c<0$,
    \begin{itemize}[leftmargin=0.5cm]
    \item $|c|\leq k+1$, and $n > 0$. Then $Q$ has \textbf{$n$-class $S_{-}$}.
    \item $|c|\leq k+2$, and $n \leq 0$. Then $Q$ has \textbf{$n$-class $S_{-}$}.
    \end{itemize}
\end{itemize}
and define $n$-classes and $n$-types inductively following Definition \ref{type_induct}. We define $\tilde{\S}$ to be the same as $\S$ for polynomials with positive leading coefficients. For polynomials with negative leading coefficients, the previous extension of Definition \ref{type_base} naturally yields an extension of Definition \ref{type_induct} and \ref{def_succ} with the following exceptions:
\begin{enumerate}
    \setcounter{enumi}{5}
    \item If $P$ is of $n$-class $-E_2$, we let $Q$ be the rewriting of $P+1$ given by Proposition \ref{rewriting}, we then set  $\tilde{\S}(P) = Q-X^{-1}$. When $P$ is of $n$-class $-E_1$, we then define $\tilde{\S}(P) = Q-1-X^{-1}$
	\item If $P=XR$ and $R$ is of $(n-1)$-class $-U^t_{-5}$, then $\tilde{\S}(P)$ is the rewriting of $P+1$ given by Proposition \ref{rewriting}.
\end{enumerate}
We do this to handle the potential change violations by \textbf{rule 4} and the sign change violations by \textbf{rule 5} in Proposition \ref{classification} respectively.

Denote the set of polynomials with a nonnegative leading coefficient $\mathcal{R}^+$. First we show that the successor function is a bijection from $\mathcal{R}^+ \cup \{0\}$ onto $\mathcal{R}^+$.

\begin{proposition}\label{succ_inj}
The successor function is a bijection from $\mathcal{R}^+\cup \{0\}$ onto $\mathcal{R}^+$.
\end{proposition}
\begin{proof}
We show this by showing there is an inverse function $\S^{-1}(P):\mathcal{R}^+\rightarrow \mathcal{R}^+\cup \{0\}$. Let $$\S^{-1}(P)=-\tilde{\S}(-P).$$
where $\tilde{\S}$ is the generalized successor function. $\tilde{\S}$ will be suppressed to $\S$ unless needed. We see that this function is well-defined. We first show that $\S^{-1}(\S(P))=P$ for each case from Definition \ref{def_succ}.

Suppose that $\S$ is a regular successor for $P$. In order for the identity to hold, we need $-(\S(-P-1))=P$, or equivalently, $\S(-P-1)=-P$. When $P$ is of $n$-initial, $n$-interior type, since $-P-1$ is of $n$-negative type, the identity holds. %The same is true when $P$ is of $n$-boundary type (P) and $n$-class $U_{-3}^t$.

Now suppose that $P$ is of $n$-negative type. We show that $-P-1$ is always of $n$-interior type. We begin with the base case:
\begin{center}
    \begin{tabular}{c|c|c|c|c}\hline
        $n$-class of $P$ & $-P-1$ & $n$-class of $-P-1$ & $n$-class of $-P$ & $n$ \\\hline
        $E_1$ & $-X+k-2$ & $S_+$ & $U_0$ &$n\leq 1$ \\
        $E_3$ & $-X+k-3$ & $S_+$ & $S_+$ &$n\leq 1$ \\
        $S_{-}$ & $0<c\leq k$ & $S_0$, $S_+$ or $U_{0}$ & $S_+$, $U_0$ or $U_{-2}$ &$n> 0$ \\
        $S_{-}$ & $0<c\leq k+1$ & $S_0$, $S_+$ or $U_{-3}^t$ & $S_+$, $U_{-3}^t$ or $U_{-5}^t$ &$n\leq 0$ \\
        \hline
    \end{tabular}.
\end{center}
Thus, the claim is true for the base case. Now, in general, if $P=X R+c$, then
\begin{center}
    \begin{tabular}{c|c|c|c|c}\hline
        $n$-class of $P$ & $(n-1)$-class of $R$& $-R$ & $c'=-c-1$ & $-P-1$\\\hline
        $E_2$ & $E_1$ & $U_0$ & $k-2$ & $S_+$\\
        $E_3$ & $E_2$ & $S_+$ & $k-2$ & $S_+$\\
        $S_-$ & $E_3$ & $S_+$ & $k-2$ & $S_+$\\
        $E_1$ & $E_2$ & $S_+$ & $k-1$ & $S_+$\\
        $E_2$ & $E_3$ & $S_+$ & $k-1$ & $S_+$\\
        $S_-$ & $S_0$ & $S_0$ & $0\leq c'\leq k$ & $S_0$, $S_+$, $U_0$\\
        $S_-$ & $*$ & $S_-$ & $0\leq c'\leq k-1$ & $S_0$, $S_+$, $U_0$\\
        \hline
    \end{tabular}.
\end{center}
Working out each case directly shows that $\S(-P-1)=-P$ as desired except when $P$ is of $n$-class $S_-$ and $R$ is of $(n-1)$-class $S_-$. We note that since $P$ is $n$-reduced, $-R$ being $n$-class $U_{-2}$ can be ruled out. Furthermore, $-P$ needs to be $n$-reduced, significantly reducing the number of cases to be considered. In this case, working out explicitly we see that
\begin{center}
    \begin{tabular}{c|c|c|c}\hline
        class of $P$ & $-R$ & $c'=-c-1$ & $-P-1$\\\hline
        $S_-$ & $S_+$ & $0\leq c'\leq k$ & $S_+$, $U_0$\\
        $S_-$ & $U_0$ & $0\leq c'\leq k-1$ & $S_+$\\
        $S_-$ & $U_{-1}$ & $0\leq c'\leq k-1$ & $S_+$, $U_0$\\
        $S_-$ & $U_{-3}^t$ & $0\leq c'\leq k$ & $S_+$\\
        $S_-$ & $U_{-4}^t$ & $0\leq c'\leq k-1$ & $S_+$, $U_{-3}^t$\\
        $S_-$ & $U_{-5}^t$ & $0\leq c'\leq k-2$ & $S_+$\\
        \hline
    \end{tabular}..
\end{center}
Hence, the successor for $-P-1$ is regular.

Now suppose that $\S$ is irregular. For the base case,  the identity holds trivially by construction.
Let $P$ be of $n$-class $U^t_{-4}$ and $U^t_{-5}$. It is easy to check that $\S(P)$ is of $n$-class $E_2$ and $E_1$, respectively. For $-\S(P)$, by construction of the generalized successor, $\tilde{\S}(-\S(P))=P$ as required.

Now let $P$ be of $n$-class $U$. For the base case, we check when $P=XR+c_0$ is of $n$-class $U_{-2}$ where $R$ is of $n$-class $S$ or $E$. For the recursive case, there exists smallest $d>2$ such that
$$P=X^d B+XA+c_0$$
where $B$ has $(n-d)$-class $S$ or $E$ and $A$ is a degree $d-2$ polynomial. Denote $P=(B|A|c_0)$. By choice of $d$, $XA+c_0+1$ corresponds to the rightmost minimal poison subword of $P+1$.
\begin{center}
    \begin{tabular}{c|c|c|c|c}\hline
        $n$-class of $P$ & $\Pot(A)$ & $c_0$ & $-c_0'$ & $n$-class of $-\S(P)$\\\hline
		$U_{-2}$ & $-2-\sign^+(c_j\cdot c_{j+1})$ & $k$ & $k$ & $U_{-2}$\\
		$U_{-1}$ & $-1-\sign^+(c_j\cdot c_{j+1})$ & $k$ & $k$ & $U_{-1}$\\
		$U_{-1}$ & $-1-\sign^+(c_j\cdot c_{j+1})$ & $k-1$ & $k-1$ & $U_{-1}$\\
        \hline
    \end{tabular}
\end{center}
We observe that $-\S(P)$ can only have $n$-class $U$ and that the rightmost minimal poison subword for $-\S(P)$ can only be on the same position as $A$. The rewriting rule for $-\S(P)$ is identical to the rewriting rule for $P$, and thus, the identity holds.

For $n$-boundary (S), let $Q=XP$ where $P$ is of $(n-1)$-class of $E_1$. $\S(Q)$ is also of the form $XR$ where $R$ is of $(n-1)$-class $-U^t_{-5}$, and the identity follows from construction of the generalized successor. When $Q=XP$ and $P$ is of $(n-1)$-class $-U_{-2}$, one can observe $-\S(Q)$ is also $n$-boundary (S) which completes the proof.

The other direction $\S(\S^{-1}(P))=\S(-\S(-P))$ follows from symmetry.

\end{proof}

\begin{proposition}\label{succ-rep}
Suppose P is a $n$-reduced polynomial with the nonnegative leading coefficient that represents $x$.
\begin{enumerate}
\item If $P$ has $n$-initial type, then $L_n(\S(P)) = L_n(P) + 1$ and $\S(P)$ represents $x+b$.
\item If $P$ has $n$-interior type, then $L_n(\S(P)) = L_n(P) + 1$ and $\S(P)$ represents $x+b$. 
\item If $P$ has $n$-negative type, then $L_n(\S(P)) = L_n(P) - 1$ and $\S(P)$ represents $x+ b$.
\item If $P \in U_{-5}^t$ or $P \in U_{-2}$, then $L_n(\S(P)) = L_n(P)$ and $\S(P)$ represents $x -a$. %x-a for ordinary successor
\item If $P \in U_{-4}^t$ or $P \in U_{-1}$, then $L_n(\S(P)) = L_n(P)$ and $\S(P)$ represents $x +b - a$.
\item If $P$ has $n$-boundary type (S), then $L_n(\S(P)) = L_n(P)$ and $\S(P)$ represents $x+b$.
\end{enumerate}
\end{proposition}

\begin{proof}
We will only prove the case when $P\in U_{-2}$. The other cases can be proved similarly. 

Suppose $P \in U_{-2}$. We then have $P = (\cdots | A)$ for some $A$ with degree $d$ and potential $-2-\sign^+(c_j\cdot c_{j+1})$. The rewriting of $P+1$ is by subtracting $(X^{d}+X^{d-1} + \cdots +1)(X-(2k+1)+X^{-1})$, and thus, we have
\begin{align*}
&\S(P) \\
=&P+1 + (X^{d}+X^{d-1} + \cdots +1)(X-(2k+1)+X^{-1}) -1-X^{-1}\\
=& P+1 + (X^{d+1} -2kX^{d}-(2k-1)X^{d-1} - \cdots -(2k-1)X -2k +X^{-1})-1-X^{-1}\\
=& P + (X^{d+1} -2kX^{d}-(2k-1)X^{d-1} - \cdots -(2k-1)X -2k)\\
=& P + X^{d+1} -X^d -1 -(2k-1)(X^{d}+X^{d-1} + \cdots +X +1).
\end{align*}

Recall that the potential of a coefficient is exactly the length change when adding $2k-1$. Therefore, the length difference in the last $d+1$ digits of $P$ and $\S(P)$ is $-2-\sign^+(c_j\cdot c_{j+1})+2 = -\sign^+(c_j\cdot c_{j+1})$, which with $c_{j+1} > 0$, is exactly the opposite of the effect of $X^{d+1}$ has on the length change. Thus, $L_n(\S(P)) = L_n(P)$. We also have $\S(P)$ represents $x-a$ which follows from the fact that rewriting does not change the element the polynomial represents and that $X^{-1}$ represents $a$. \end{proof}

%\begin{proposition}
%Suppose P is a $n$-reduced polynomial with the positive leading coefficient that represents $x$.
%\begin{enumerate}
%\item If $P$ has $n$-initial $n$-type, then $L_n(\S(P)) = L_n(P) + 1$ and $\S(P)$ represents $x+b$.
%\item If $P \in S_0$ or $P \in S_+$, then $L_n(\S(P)) = L_n(P) + 1$ and $\S(P)$ represents $x+b$. 
%\item If $P \in S_-$, then $L_n(\S(P)) = L_n(P) - 1$ and $\S(P)$ represents $x+ b$.
%\item If $P \in U_{-3}^t$ or $ P \in U_{0}$, then $L_n(\S(P)) = L_n(P) + 1$ and $\S(P)$ represents $x + b$.
%\item If $P \in U_{-4}^t$ or $P \in U_{-1}$, then $L_n(\S(P)) = L_n(P)$ and $\S(P)$ represents $x +b - a$.
%\item If $P \in U_{-5}^t$ or $P \in U_{-2}$, then $L_n(\S(P)) = L_n(P)$ and $\S(P)$ represents $x+b -a$. %x-a for ordinary successor
%\item If $P$ is $n$-class $E_\bullet$, then $L_n(\S(P)) = L_n(P)$ and $\S(P)$ represents $x +b$.
%\end{enumerate}
%\end{proposition}

For a Laurent polynomial $F(X)$, we will write $\overline{F}(X) = X^{-1}\cdot F(X^{-1})$. We define the successor function on its principal part $Q(X)$ dually to the polynomial part, i.e., the $n$-successor of $Q(X)$ is $\overline{\S(\overline{Q}(X^{-1})}$, where $\S$ is the $(-n)$-successor function.

Our main interest is in quantifying group elements, not $n$-reduced polynomials. Suppose $F(X) = \sum_{j=-p}^{q} c_j X^j$ represents $x$, and write its polynomials and principal parts as $P(X)$ and $Q(X)$, respectively. Then since the successor function is surjective, we can write $P(X) = \pm\S^n(0)$ and $\overline{Q}(X) = \pm\S^m(0)$. Thus, $(\pm\S^n(0)|\pm\overline{\S^m(0)})$ represents $x$. However, for a fixed $n$, there may be more than one way to write it as $(\pm\S^n(0)|\pm\overline{\S^m(0)})$. To quantify this, we will need the \emph{proof} of the following lemma even though the lemma is not used directly. 

We view the sequence $\{\S^i(0)\}$ of successive polynomials as a sequence of points $g_i$ they represent in $\mathbb{Z}^2=\left< a,b \right>$, and the successor corresponds to moving in either $b$, $-a$ or $b-a$ direction depending on the type. Note that if $\overline{Q}(X)$ represents $ka+tb$ then $Q(X)$ represents $ta+kb$. We state the lemma in a slightly more general setting.

\begin{lemma}[Spanning lemma]\label{Exi.Rep}
Let $G = \left< a,b \right>$ be a free abelian group. Suppose there exists an infinite sequence $\{g_i \} \subset G$ such that $g_0 = 0$, $g_1 = b$ and $g_{i+1} \in \{g_{i} + b, g_i + b -a, g_i-a\}$. Moreover, assume that if $g_{i+1} - g_i \in \{b - a,-a\}$, then $g_{i+2} - g_{i+1} = g_{i+3}-g_{i+2}=  b$. We also define $\{h_i\}$ such that if $g_i = ka+ tb$, then $h_i = ta+ kb$. Then for all $\gamma \in G$, there exist $i,j \in \mathbb{N}$ such that $\gamma = \pm g_{i} \pm h_{j}$.
\end{lemma}

\begin{proof}
We first prove that if some $\gamma$ is a sum of some $g_i$ and $h_j$, then any $\gamma + ma+nb$ is also for any $m,n\in \N$. We proceed by induction on $n$ and $m$. Suppose that $\gamma = g_i+h_j$. We first need to find a pair $(i',j')$ such that $\gamma + b = g_{i'}+h_{j'}$. If $g_{i+1}=g_i+b$, we are done. If $g_{i+1}=g_i+(b-a)$, then
    \begin{align*}
        g_{i+1}+h_{j+1}&=g_i+(b-a)+h_j+(a),\\ g_{i+2}+h_{j+2}&=g_i+(b-a)+(b)+h_j+(-b)+(a) \text{, or}\\
        g_{i+2}+h_{j+2}&=g_i+(b-a)+(b)+h_j+(a-b).
    \end{align*}
On the other hand, if $g_{i+1}=g_i+(-a)$, then
    \begin{align*}
        g_{i+2}+h_{j+1}&=g_i+(-a)+(b)+h_j+(a),\\ g_{i+2}+h_{j+2}&=g_i+(b-a)+(b)+h_j+(-b)+(a) \text{, or}\\ g_{i+3}+h_{j+2}&=g_i+(-a)+(b)+(b)+h_j+(-b)+(a).
    \end{align*}
We observe that by the symmetry between $g_i$ and $h_i$, the $a$'s can be incremented in the same manner. This proves the lemma for elements in $\{g_i+ma\}\cup\{h_j+nb\}$, i.e., the upper-right region cut out by the trajectories of $g_i$ and $h_i$ (see figure \ref{eigne}). By switching signs, any element in $\{-g_i-ma\}\cup\{-h_j-nb\}$ can be expressed as $-g_i-h_j$ as well.

\begin{figure}[ht!]\label{eigne}
	\includegraphics[scale=0.6]{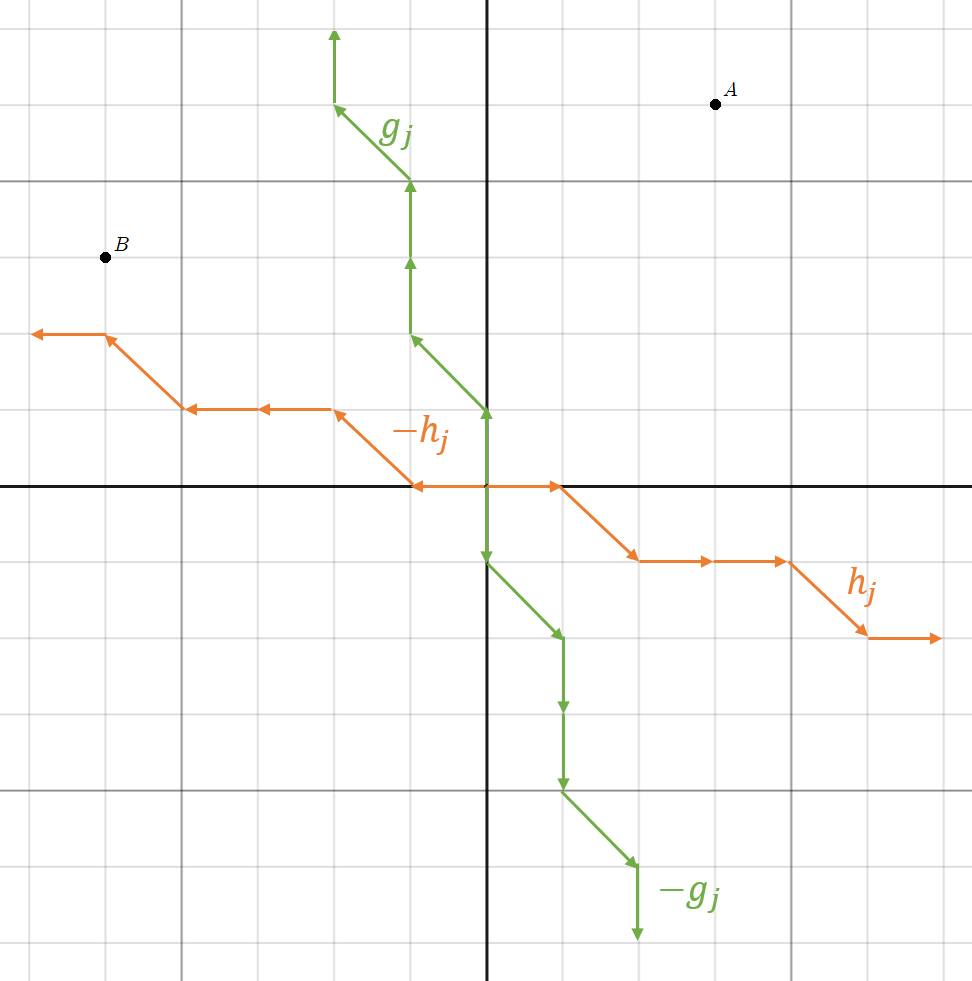}
	\caption{For point $A$, there is some $(i,j)$ such that $g_i+h_j$ represents $A$. For point $B$, the representative to look for is $g_i-h_j$. Instead of actually finding $(i,j)$, we increment by $\pm a$ and $\pm b$ from a known pair.}
\label{eigen}
\end{figure}

By symmetry, this leaves out the case when $\gamma\in \{g_i-ma\}\cap\{-h_j+nb\}$. We will show that such $\gamma$ can be written as $g_i-h_j$ for some $i,j\in N$. By induction, let $\gamma=g_i-h_j$ for $j>1$. We need to find a pair $(i',j')$ such that $\gamma+b=g_{i'}-h_{j'}$ where $i'>i$ and $j'<j$. 
%When $i=0$ and $\gamma=-h_j=ma+nb$, then $ma+(n+1)b=g_1-h_j$ so the claim can be seen to be true. We work when $i=1$. Suppose $\gamma=g_1-h_j=ma+nb$. 
If $g_{i+1}-g_i=b$, then $\gamma+b=g_i-h_j$, and we are done. If $g_{i+1}-g_i=b-a$, then either
    \begin{align*}
        g_{i+1}-h_{j-1}&=g_i+(b-a)-h_j+(a),\\
        g_{i+2}-h_{j-1}&=g_i+(b-a)+(b)-h_j+(a-b) \text{, or}\\
        g_{i+2}-h_{j-2}&=g_i+(b-a)+(b)-h_j+(-b)+(a) \quad (j>1)
    \end{align*}
where the last case only appears when $j\geq 2$ as $h_1 = a$. On the other hand, if $g_{i+1}-g_i=-a$, then either
    \begin{align*}
        g_{i+1}-h_{j-1}&=g_i+(-a)+(b)-h_j+(a),\\ 
	g_{i+3}-h_{j-1}&=g_i+(-a)+(b)+(b)-h_j+(a-b)\text{, or}\\
        g_{i+3}-h_{j-2}&=g_i+(-a)+(b)+(b)-h_j+(-b)+(a). \quad (j>1)
    \end{align*}
%In fact, for any $g_i-h_j$, this can be done repeatedly until it lands on the orbit $\{g_{k}\}$. 
When $j = 0$, $\gamma = g_i$ and $\gamma+b$ will either be $g_{i+1}$ or is covered by the first case of the proof. Hence, any element in the set $\{g_i-ma\}\cap\{-h_j+nb\}$ can be expressed as a difference $g_i-h_j$. By symmetry, this can be done for the set $\{-g_j+ma\}\cap \{h_j-nb\}$.
\end{proof}

For counting, we need to know how many representatives of a single $\gamma$ exist. We use the same framework as in the Lemma \ref{Exi.Rep} to figure out all the cases and then find the appropriate types and classes in relation to the successor function.

As in the lemma, we first look at the sum of two sequences $g_i$ and $h_j$. We find all possible pairs of sequences satisfying $$g_i+h_j=g_{i'}+h_{j'}.$$ Equivalently, we can find sequences with weaker conditions whose sum is 0. Here we only require $-a$ and $b-a$ to be followed by two consecutive $b$'s in the sequence $\{g_i\}$ and a similar requirement for the sequence $\{h_j\}$. When any of two sequences is long, the $b$'s (or $a$'s) dominate the sequence. Hence, it suffices to look at combinations shorter than 3 symbols:
\begin{align*}
    0=& (-a)+(a)=(b)+(-b)=(b-a)+(a-b)\\
    =&(b-a)+(a+(-b))=(b-a)+((-b)+a)\\
    =&(b+(-a))+(a+(-b))=((-a)+b)+((-b)+a)=(b+(-a))+((-b)+a)\\
    =&((-a)+b)+(a+(-b)).    
\end{align*}

Now we consider the Laurent polynomial $F(X) = \sum_{j=-p}^{q} c_j X^j$ that represents $\gamma$, and write its polynomials and principal parts as $P(X)$ and $Q(X)$, respectively. Suppose both $P$ and $\overline{Q}$ have nonnegative leading coefficients. For $n\geq 0$, when $P$ is an $n$-reduced polynomial and $\overline{Q}$ is a $(-n)$-reduced polynomial, we have $P$ representing some $g_i$ and $Q$ representing some $h_j$ and $x = g_i+h_j$. Thus, the overcountings are in the form $g_i + h_j = g_{i'}+h_{j'}$, and based on the discussion above and Prop \ref{succ-rep} can now be classified in terms of types and classes.

$(P|\overline{Q})=(\S(P)|\overline{\S(Q)})$ occurs when

\vspace{1ex}
\begin{center}
    \begin{tabular}{c|c}\hline
        $n$-classes of $P$ & $(-n)$-classes of $\overline{Q}$\\\hline
        $U_{-2}$, $U^t_{-5}$ & $0$, $S_{\bullet}$, $U_{0}$, $E_{\bullet}$, $U^t_{-3}$ \\\hline
        $0$, $S_{\bullet}$, $E_{\bullet}$, $U_{0}$, $U^t_{-3}$ & $U_{-2}$, $U^t_{-5}$ \\\hline
        $U_{-1}$, $U^t_{-4}$ & $U_{-1}$, $U^t_{-4}$ \\\hline
    \end{tabular}.
\end{center}
Similarly, $(P|\overline{Q})=(\S(P)|\overline{\S^2(Q)})$ occurs when
\vspace{1ex}
\begin{center}
    \begin{tabular}{c|c|c}\hline
        $n$-classes of $P$ & $(-n)$-classes of $\overline{Q}$ & $(-n)$-classes of $\overline{\S(Q)}$\\\hline
        $U_{-1}$, $U^t_{-4}$ & $U_{0}$, $U^t_{-3}$ & $U_{-2}$, $U^t_{-5}$ \\\hline
        $U_{-1}$, $U^t_{-4}$ & $U_{-2}$, $U^t_{-5}$ & $S_{-}$, $E_{\bullet}$ \\\hline
    \end{tabular}.
\end{center}
And finally, $(P|\overline{Q})=(\S^2(P)|\overline{\S^2(Q)})$ can be classified.

\vspace{1ex}
\begin{center}
    \begin{tabular}{c|c|c|c}\hline
        $P$ & $\S(P)$ & $\overline{Q}$ & $\overline{\S(Q)}$ \\\hline
        $U_{0}$, $U^t_{-3}$ & $U_{-2}$, $U^t_{-5}$ & $U_{0}$, $U^t_{-3}$ & $U_{-2}$, $U^t_{-5}$ \\\hline
        $U_{-2}$, $U^t_{-5}$ & $S_{-}$, $E_{\bullet}$ & $U_{-2}$, $U^t_{-5}$ & $S_{-}$, $E_{\bullet}$ \\\hline
        $U_{0}$, $U^t_{-3}$ & $U_{-2}$, $U^t_{-5}$ & $U_{-2}$, $U^t_{-5}$ & $U_{0}$, $U^t_{-3}$ \\\hline
        $U_{-2}$, $U^t_{-5}$ & $S_{-}$, $E_{\bullet}$ & $U_{0}$, $U^t_{-3}$ & $U_{-2}$, $U^t_{-5}$ \\\hline
    \end{tabular}.
\end{center}

Now we consider the case when $P$ is positive and $\overline{Q}$ is negative, so $P+Q$ represents $x = g_i-h_j$:
\begin{align*}
    b&=(b)-0=0-(-b)\\
    -a&= (-a)-0=0-(a)\\
    b-a&=(b-a)-0=0-(a-b)    .
\end{align*}
In terms of the successor, this corresponds to $$(\S(P)|\overline{Q})=(P|\overline{\S(Q)}),$$
and the possible classes for $P$ and $Q$ are described below.
\vspace{1ex}
\begin{center}
    \begin{tabular}{c|c}\hline
        $n$-classes of $P$ & $(-n)$-classes of $\overline{Q}$\\\hline
        $0$, $S_{\bullet}$, $U_{0}$, $E_{\bullet}$, $U^t_{-3}$ & $0$, $S_{\bullet}$, $U_{0}$, $E_{\bullet}$, $U^t_{-3}$\\\hline
        $U_{-2}$, $U^t_{-5}$ & $U_{-2}$, $U^t_{-5}$ \\\hline
        $U_{-1}$, $U^t_{-4}$ & $U_{-1}$, $U^t_{-4}$\\\hline
    \end{tabular}
\end{center}

\section{Growth series for $n$-reduced polynomials}
Recall that the $n$-length $L_n$ of a Laurent polynomial representation $$P(X)=\sum_{j=-p-1}^{q} c_j X^j$$ of an element $g\in G$ is defined as
$$
L_n(P) = \|g\|=2p + 2q - |n| + \sum_{j=-p-1}^q |c_j|
$$
where $p \geq \text{max}\{0,-n\}$ and $q \geq \text{max}\{0,n\}$. Here we assumed the polynomial had zero as its coefficients on both ends if these two conditions are not met. In particular, for $n$-reduced polynomials, we can restrict to
$$
L_n(P) =2q-n + \sum_{j=0}^q |c_j|
$$
and pass $-|n|$ to either the polynomial part or the principal part.

We define the generating function for polynomials of given $n$-class of degree $d$ with a nonnegative leading coefficient and use the inductive definition to find the recursive relation. Let
$$S_{+}^{n,d}(t)=\sum_{P\in S_+,\:deg(P)= d} t^{L_n(P)}$$
be the generating function for the polynomials of $n$-class $S_{+}$ of degree $d$, and for all stable classes, denote
$$S^{n,d}(t)=\begin{pmatrix}S_{+}^{n,d}(t)\\S_{0}^{n,d}(t)\\S_{-}^{n,d}(t)\end{pmatrix}.$$
In a similar fashion, let
\begin{align*}
    U^{n,d}(t)=\begin{pmatrix}U_{0}^{n,d}(t)\\U_{-1}^{n,d}(t)\\U_{-2}^{n,d}(t)\end{pmatrix}, \quad T^{n,d}(t)=\begin{pmatrix}T_{-3}^{n,d}(t)\\T_{-4}^{n,d}(t)\\T_{-5}^{n,d}(t)\end{pmatrix}, \quad E^{n,d}(t)=\begin{pmatrix}E_{1}^{n,d}(t)\\E_{2}^{n,d}(t)\\E_{3}^{n,d}(t)\end{pmatrix}
\end{align*}
be the generating functions for $n$-class $U$, $n$-class $U^t$ and $n$-class $E$, respectively. We observe that if $n\leq 0$, then
\begin{align*}
    S^{n,d}(t)=t^{-n} S^{0,d}(t), \quad U^{n,d}(t)=t^{-n} U^{0,d}(t)\\
    T^{n,d}(t)=t^{-n} T^{0,d}(t), \quad E^{n,d}(t)=t^{-n} E^{0,d}(t)
\end{align*}
since the rules for $n$-reduced polynomials are identical when $n\leq 0$ but the length changes by $|n|$. We will be using this identity repeatedly to obtain the generating functions when $n=0$ without the degree constraint.

We begin computing these vectors with small $d$ listed in Definition \ref{type_base}.
Let $d=0$. When $n>0$, we have $n$-classes $S_0$, $S_+$, $U_0$ and $U_{-2}$, and $L_n(c_0)=n+|c_0|$. On the other hand, if $n\leq 0$, we have $n$-classes $S_+$, $U^t_{-3}$ and $U^t_{-5}$. We do not include the initial type since it plays no role in the recursive definition of $n$-types. In this case, $L_n(c_0)=-n+|c_0|$. So for stable $n$-types, we have

\begin{align*}
    S^{0,0}(t)&=\begin{pmatrix}t+\cdots+t^{k} \\0 \\0\end{pmatrix}=\frac{1}{t-1}\begin{pmatrix}t^k-t\\0\\0\end{pmatrix}\\
    S^{n,0}(t)&=t^{n}\begin{pmatrix}t+\cdots+t^{k-1} \\1 \\0\end{pmatrix}=\frac{1}{t-1}\begin{pmatrix}t^{k+1}-t\\t-1\\0\end{pmatrix}
    \text{ if $n>0.$}
\end{align*}

Similarly,
\begin{align*}
    U^{0,0}(t)=0, \quad T^{0,0}(t)=\begin{pmatrix}t^{k+1} \\0\\t^{k+2}\end{pmatrix}\\
    U^{n,0}(t)=\begin{pmatrix}t^k \\0 \\t^{k+1}\end{pmatrix}, \quad T^{n,0}(t)=0\text{  if $n>0$.}\\
    E^{n,0}(t)=0\text{  for all $n$.}
\end{align*}

The $n$-class $E$ appears only when the degree is greater than 1 and $n< d$. Since $P(X)=X-k+1$ is $0$-class $E_1$ and $P(X)=X-k+2$ is $0$-class $E_3$, we have that
\begin{align*}
    E^{n,1}(t)&=0 \text{ if $n>0$}\\
    E^{0,1}(t)&=\begin{pmatrix}t^{k+2} \\0 \\t^{k+1}\end{pmatrix}.
\end{align*}

Now we introduce recursion using Definition \ref{type_induct}. We begin with $n$-class $E$.
\begin{align*}
    E^{(n+1,d+1)}(t)&=t
        \begin{bmatrix}
            0 & t^{k}& 0\\
            t^{k-1} & 0 & t^{k} \\
            0 & t^{k-1} & 0
        \end{bmatrix}
    E^{n,d}(t)=P_{E,E}E^{n,d}.
\end{align*}
Similarly for other classes except stable classes, we have
\begin{align*}
    U^{(n+1,d+1)}(t)&=
        t\begin{bmatrix}
            t^k&t^k&t^{k-1}\\
            0&0&0\\
            t^{k+1}&t^{k+1}&t^{k+1}\\
        \end{bmatrix}
    S^{n,d}(t)
    +
        t\begin{bmatrix}
            0&t^{k-1}&0\\
            t^{k}&0&t^{k-1}\\
            0&t^k&0\\
        \end{bmatrix}
    U^{n,d}\\
    &=P_{U,S}S^{n,d}+P_{U,U}U^{n,d}\\
    T^{(n+1,d+1)}&=
        t\begin{bmatrix}
            0&t^{k-1}&0\\\
            t^{k}&0&t^{k-1}\\
            0&t^{k}&0
        \end{bmatrix}
    T^{n,d}=P_{T,T}T^{n,d}.
\end{align*}

For stable classes, the same argument can be done, giving us
\begin{align*}
    S^{(n+1,d+1)}&=\frac{t}{t-1}
        \begin{bmatrix}
            t^k-t&t^k-t&t^{k-1}-t\\
            t-1&t-1&t-1\\
            * & * & *\\
        \end{bmatrix}
    S^{n,d}\\
    &+\frac{t}{t-1}
        \begin{bmatrix}
            t^{k}-t&t^{k-1}-t&t^{k-1}-t\\
            t-1&t-1&t-1\\
            * & * & *\\
        \end{bmatrix}
    U^{n,d}\\
    &+\frac{t}{t-1}
        \begin{bmatrix}
            t^{k}-t&t^{k-1}-t&t^{k-1}-t\\
            t-1&t-1&t-1\\
            * & * & *\\
        \end{bmatrix}
    T^{n,d}\\
    &+\frac{t}{t-1}
        \begin{bmatrix}
            t^{k}-t&t^{k}-t&t^{k}-t\\
            0&t-1&t-1\\
            t^{k-1}-t&t^{k-1}-t&t^{k}-t\\
        \end{bmatrix}
    E^{n,d}\\
    &=P_{S,S}S^{n,d}+P_{S,U}U^{n,d}+P_{S,T}T^{n,d}+P_{S,E}E^{n,d}.
\end{align*}
The recursive definition does not explicitly state what coefficients can be attached at the end to get $n$-class $S_0$. Instead of pursuing the exact condition, we note that for any polynomial $P$ of $n$-class $E$ or $S_{-}$, we know exactly the $n$-class of $-P$ by Proposition \ref{succ_inj}. Exploiting this symmetry, we have
$$S_-^{(n,d)}+\sum_j E_j^{(n,d)}=\begin{cases}S_+^{(n,d)}+S_0^{(n,d)}+U_0^{(n,d)}+U_{-2}^{(n,d)} \quad \text{if $n\geq d$}\\
												S_+^{(n,d)}+S_0^{(n,d)}+T_{-3}^{(n,d)}+T_{-5}^{(n,d)}\quad \text{if $n<d$}\end{cases}$$

Denoting
\begin{align*}
    \mathcal{P}&=
    \begin{bmatrix}
        P_{S,S}&P_{S,U}&P_{S,T}&P_{S,E}\\
        P_{U,S}&P_{U,U}&0&0\\
        0&0&P_{T,T}&0\\
        0&0&0&P_{E,E}\\
    \end{bmatrix}.
\end{align*}
by the mentioned symmetry, we have the following identity

\begin{align*}
	\mathcal{P} \begin{pmatrix}S^{(n,d)}\\U^{(n,d)}\\T^{(n,d)}\\E^{(n,d)}\\\end{pmatrix}
	 = \begin{bmatrix}   A & B & C & D\\
						P_{U,S}&P_{U,U}&0&0\\
        					0&0&P_{T,T}&0\\
        					0&0&0&P_{E,E}\\\end{bmatrix} \begin{pmatrix}S^{(n,d)}\\U^{(n,d)}\\T^{(n,d)}\\E^{(n,d)}\\\end{pmatrix}
\end{align*}

where
\begin{align*}
A&= \frac{t}{t-1}
        \begin{bmatrix}
            t^k-t&t^k-t&t^{k-1}-t\\
            t-1&t-1&t-1\\
            t^k-1 & t^k-1 & t^{k-1}-1\\
        \end{bmatrix}+P_{U,S}\\
B&= \frac{t}{t-1}
        \begin{bmatrix}
            t^{k}-t&t^{k-1}-t&t^{k-1}-t\\
            t-1&t-1&t-1\\
            t^{k}-1&t^{k-1}-1&t^{k-1}-1\\
	  \end{bmatrix}+\begin{bmatrix}
            1&0&0\\
            0&0&0\\
            0&0&1\\
        \end{bmatrix}P_{U,U}\\
C&=\frac{t}{t-1}
        \begin{bmatrix}
            t^{k}-t&t^{k-1}-t&t^{k-1}-t\\
            t-1&t-1&t-1\\
            t^{k}-1&t^{k-1}-1&t^{k-1}-1\\
        \end{bmatrix}+\begin{bmatrix}
            1&0&0\\
            0&0&0\\
            0&0&1\\
        \end{bmatrix}P_{T,T}\\
D&=\frac{t}{t-1}
        \begin{bmatrix}
            t^{k}-t&t^{k}-t&t^{k}-t\\
            0&t-1&t-1\\
            t^{k-1}-t&t^{k-1}-t&t^{k}-t\\
        \end{bmatrix}-P_{E,E}
\end{align*}
amd we will use this new block matrix instead of $\mathcal{P}$.

We sum over all possible degrees
\begin{align*}
    \begin{pmatrix}
        S^{0}\\
        U^{0}\\
        T^{0}\\
        E^{0}\\
    \end{pmatrix}=
		\begin{pmatrix}
                S^{0,1}\\
                U^{0,1}\\
                T^{0,1}\\
                E^{0,1}\\
            \end{pmatrix}+\sum_{d=0}^\infty \mathcal{P}^d
            \begin{pmatrix}
                S^{0,1}\\
                U^{0,1}\\
                T^{0,1}\\
                E^{0,1}\\
            \end{pmatrix}=
		   \begin{pmatrix}
                S^{0,1}\\
                U^{0,1}\\
                T^{0,1}\\
                E^{0,1}\\
            \end{pmatrix}+
		   (1-\mathcal{P})^{-1}
            \begin{pmatrix}
                S^{0,1}\\
                U^{0,1}\\
                T^{0,1}\\
                E^{0,1}\\
            \end{pmatrix},
\end{align*}
and for $n$-classes, we have
\begin{align*}
    \begin{pmatrix}
        S^{n}\\
        U^{n}\\
        T^{n}\\
        E^{n}\\
    \end{pmatrix}=(\mathcal{P}/t)^n
    \begin{pmatrix}
        S^{0}\\
        U^{0}\\
        T^{0}\\
        E^{0}\\
    \end{pmatrix}.
\end{align*}
Since $\mathcal{P}$ is a $12\times12$ matrix with polynomial entries, $(1-\mathcal{P})^{-1}$ is a matrix with rational function entries. Thus, $S^n$, $U^n$, $T^n$ and $E^n$ are all rational for fixed $n$. $(1-\mathcal{P}/t)^{-1}$ is also a matrix with rational function entries, hence the sum of $S^n$, $U^n$, $T^n$ and $E^n$ over all $n$ are rational. Since every \textit{$n$-minimal} polynomial representing a group element can be decomposed into a $n$-reduced polynomial and $-n$-reduced polynomial with possibly multiplicities described in Proposition \ref{Exi.Rep}, the set of all $n$-minimal polynomials exhibits a rational growth. Furthermore, taking account of the multiplicities in representing group elements only changes the end function by rational functions because $U^n$, $T^n$ and $E^n$ are all rational.

%\section{Possible generalizations}

%Do this for groups where $t$ acts on $\Z^{2n}$ via a matrix with $2\times 2$ blocks.

\bibliographystyle{amsalpha}
\bibliography{ref}

\end{document}